\documentclass[11pt]{amsart}
\usepackage{latexsym,graphicx}

\numberwithin{equation}{section}
\theoremstyle{plain}

\theoremstyle{remark}

\theoremstyle{definition}

\newcommand{\D}{{\mathcal D}}
\newcommand{\E}{\mathcal E}

\newcommand{\G}{{\mathcal G}}

\newcommand{\K}{{\mathcal K}}
\renewcommand{\L}{{\mathcal L}}
\newcommand{\M}{{\mathcal M}}
\newcommand{\N}{\mathbb N}

\newcommand{\V}{{\mathcal V}}

\newcommand{\dist}{\operatorname{dist}}

\newcommand{\fp}{\operatorname{FP}}

\newcommand{\Int}{\operatorname{Int}}

\renewcommand{\span}{\operatorname{span}}

\newcommand{\supp}{\operatorname{Supp}}

\def\half{{1 \over 2}}

\newcommand{\oa}{\overrightarrow}

\newcommand{\ol}{\overline}

\def\XXint#1#2#3{{\setbox0=\hbox{$#1{#2#3}{\int}$}
      \vcenter{\hbox{$#2#3$}}\kern-.5\wd0}}

\begin{document}

\def\cal{\mathcal}

\font\tpt=cmr10 at 12 pt
\font\fpt=cmr10 at 14 pt

\font \fr = eufm10




\overfullrule=0in

\def\boxit#1{\hbox{\vrule
 \vtop{%
  \vbox{\hrule\kern 2pt %
     \hbox{\kern 2pt #1\kern 2pt}}%
   \kern 2pt \hrule }%
  \vrule}}

  \def\harr#1#2{\ \smash{\mathop{\hbox to .3in{\rightarrowfill}}\limits^{\scriptstyle#1}_{\scriptstyle#2}}\ }

\def\AA{1}
\def\BB{2}
\def\CC{3}
\def\DD{4}
\def\EE{5}
\def\FF{6}
\def\GGG{7}
\def\HH{8}
\def\II{9}
\def\JJ{10}
\def\KK{11}
\def\LL{12}
\def\MM{13}

\def\ALL{1}
\def\BTA{2}
\def\BL{3}
\def\BRE{4}
\def\CNS{5}
\def\CIL{6}
\def\CRA{7}
\def\DDD{8}
\def\DDR{9}
\def\GEO{10}
\def\HYP{11}
\def\BEL{12}
\def\AC{13}
\def\SURVEY{14}
\def\NOTES{15}
\def\AET{16}
\def\LAG{17}
\def\KRY{18}
\def\PLI{19}
\def\RT{20}
\def\SLO{21}
\def\TRU{22}
\def\TWC{23}
\def\WAL{24}

 \def\GG{{{\bf G} \!\!\!\! {\rm l}}\ }

\def\GL{{\rm GL}}

\def\bll{I \!\! L}

\def\IFF{\qquad\iff\qquad}
\def\bra#1#2{\langle #1, #2\rangle}
\def\bbf{{\bf F}}
\def\bbj{{\bf J}}
\def\Jtn{{\bbj}^2_n}  \def\JtN{{\bbj}^2_N}  \def\JoN{{\bbj}^1_N}
\def\jt{j^2}
\def\jtx{\jt_x}
\def\Jt{J^2}
\def\Jtx{\Jt_x}
\def\bpp{{\bf P}^+}
\def\bpt{{\wt{\bf P}}}
\def\fsh{$F$-subharmonic }
\def\mo{monotonicity }
\def\jet{(r,p,A)}
\def\ss{\subset}
\def\sse{\subseteq}
\def\half{\hbox{${1\over 2}$}}
\def\smfrac#1#2{\hbox{${#1\over #2}$}}
\def\oa#1{\overrightarrow #1}
\def\dim{{\rm dim}}
\def\dist{{\rm dist}}
\def\codim{{\rm codim}}
\def\deg{{\rm deg}}
\def\rank{{\rm rank}}
\def\log{{\rm log}}
\def\Hess{{\rm Hess}}
\def\Hessyp{{\rm Hess}_{\rm SYP}}
\def\trace{{\rm trace}}
\def\tr{{\rm tr}}
\def\max{{\rm max}}
\def\min{{\rm min}}
\def\span{{\rm span\,}}
\def\Hom{{\rm Hom\,}}
\def\det{{\rm det}}
\def\End{{\rm End}}
\def\Sym{{\rm Sym}^2}
\def\diag{{\rm diag}}
\def\pt{{\rm pt}}
\def\Spec{{\rm Spec}}
\def\pr{{\rm pr}}
\def\Id{{\rm Id}}
\def\Grass{{\rm Grass}}
\def\Herm#1{{\rm Herm}_{#1}(V)}
\def\arr{\longrightarrow}
\def\supp{{\rm supp}}
\def\Link{{\rm Link}}
\def\Wind{{\rm Wind}}
\def\Div{{\rm Div}}
\def\vol{{\rm vol}}
\def\foral{\qquad {\rm for\ all\ \ }}
\def\fpsh{{\cal PSH}(X,\f)}
\def\Core{{\rm Core}}
\def\dis{f_M}
\def\Re{{\rm Re}}
\def\rn{\bbr^n}
\def\pp{\cp^+}
\def\plp{\cp_+}
\def\Int{{\rm Int}}
\def\cix{C^{\infty}(X)}
\def\Gr#1{G(#1,\rn)}
\def\Symn{{\Sym(\rn)}}
\def\SymN{{\Sym(\bbr^N)}}
\def\Gpn{G(p,\rn)}
\def\fd{{\rm free-dim}}
\def\SA{{\rm SA}}
 \def\cd{{\cal C}}
 \def\cdt{{\widetilde \cd}}
 \def\cm{{\cal M}}
 \def\cmt{{\widetilde \cm}}

\def\Theorem#1{\medskip\noindent {\bf THEOREM \bf #1.}}
\def\Prop#1{\medskip\noindent {\bf Proposition #1.}}
\def\Cor#1{\medskip\noindent {\bf Corollary #1.}}
\def\Lemma#1{\medskip\noindent {\bf Lemma #1.}}
\def\Remark#1{\medskip\noindent {\bf Remark #1.}}
\def\Note#1{\medskip\noindent {\bf Note #1.}}
\def\Def#1{\medskip\noindent {\bf Definition #1.}}
\def\Claim#1{\medskip\noindent {\bf Claim #1.}}
\def\Conj#1{\medskip\noindent {\bf Conjecture \bf    #1.}}
\def\Ex#1{\medskip\noindent {\bf Example \bf    #1.}}
\def\Qu#1{\medskip\noindent {\bf Question \bf    #1.}}
\def\Exercise#1{\medskip\noindent {\bf Exercise \bf    #1.}}

\def\HoQu#1{ {\AAA T\BBB HE\ \AAA H\BBB ODGE\ \AAA Q\BBB UESTION \bf    #1.}}

\def\pf{\medskip\noindent {\bf Proof.}\ }
\def\qed{\hfill  $\vrule width5pt height5pt depth0pt$}
\def\equdef{\buildrel {\rm def} \over  =}
\def\qedqed{\hfill  $\vrule width5pt height5pt depth0pt$ $\vrule width5pt height5pt depth0pt$}
\def\mathqed{  \vrule width5pt height5pt depth0pt}

\def\V{W}

\def\df{d^{\phi}}
\def\hk{\_{\rm l}\,}
\def\n{\nabla}
\def\w{\wedge}

\def\cu{{\cal U}}   \def\cc{{\cal C}}   \def\cb{{\cal B}}  \def\cz{{\cal Z}}
\def\cv{{\cal V}}   \def\cp{{\cal P}}   \def\ca{{\cal A}}
\def\cw{{\cal W}}   \def\co{{\cal O}}
\def\ce{{\cal E}}   \def\ck{{\cal K}}
\def\ch{{\cal H}}   \def\cm{{\cal M}}
\def\cs{{\cal S}}   \def\cn{{\cal N}}
\def\cd{{\cal D}}
\def\cl{{\cal L}}
\def\cp{{\cal P}}
\def\cf{{\cal F}}
\def\ccr{{\cal  R}}

\def\gerG{{\fr{\hbox{g}}}}
\def\gerB{{\fr{\hbox{B}}}}
\def\gerR{{\fr{\hbox{R}}}}
\def\p#1{{\bf P}^{#1}}
\def\vf{\varphi}

\def\wt{\widetilde}
\def\wh{\widehat}

\def\and{\qquad {\rm and} \qquad}
\def\arr{\longrightarrow}
\def\ol{\overline}
\def\bbr{{\mathbb R}}\def\bbh{{\mathbb H}}\def\bbo{{\mathbb O}}
\def\bbc{{\mathbb C}}
\def\bbq{{\mathbb Q}}
\def\bbz{{\mathbb Z}}
\def\bbp{{\mathbb P}}
\def\bbd{{\mathbb D}}

\def\a{\alpha}
\def\b{\beta}
\def\d{\delta}
\def\e{\epsilon}
\def\f{\phi}
\def\g{\gamma}
\def\k{\kappa}
\def\l{\lambda}
\def\o{\omega}

\def\s{\sigma}
\def\x{\xi}
\def\z{\zeta}

\def\D{\Delta}
\def\L{\Lambda}
\def\G{\Gamma}
\def\O{\Omega}

\def\bd{\partial}
\def\bdf{\partial_{\f}}
\def\lag{Lagrangian}
\def\psh{plurisubharmonic }
\def\ph{pluriharmonic }
\def\pph{partially pluriharmonic }
\def\omp{$\omega$-plurisubharmonic \ }
\def\ffl{$\f$-flat}
\def\PH#1{\widehat {#1}}
\def\lloc{L^1_{\rm loc}}
\def\dbar{\ol{\partial}}
\def\lp{\Lambda_+(\f)}
\def\lpp{\Lambda^+(\f)}
\def\bo{\partial \Omega}
\def\Ob{\overline{\O}}
\def\fc{$\phi$-convex }
\def\PSH{{ \rm PSH}}
\def\SH{{\rm SH}}
\def\totr{ $\phi$-free }
\def\BM{\lambda}
\def\Der{D}
\def\CH{{\cal H}}
\def\RH{\overline{\ch}^\f }
\def\pconv{$p$-convex}
\def\MA{MA}
\def\lagpsh{Lagrangian plurisubharmonic}
\def\hermsk{{\rm Herm}_{\rm skew}}
\def\PSHl{\PSH_{\rm Lag}}
 \def\ppsh{$\pp$-plurisubharmonic}
\def\fp{$\pp$-plurisubharmonic }
\def\fh{$\pp$-pluriharmonic }
\def\Symn{\Sym(\rn)}
 \def\ci{C^{\infty}}
\def\USC{{\rm USC}}
\def\LSC{{\rm LSC}}
\def\fa{{\rm\ \  for\ all\ }}
\def\ppc{$\pp$-convex}
\def\cpt{\wt{\cp}}
\def\ft{\wt F}
\def\ob{\overline{\O}}
\def\Be{B_\e}
\def\K{{\rm K}}

\def\M{{\bf M}}
\def\N#1{C_{#1}}
\def\ds{Dirichlet set }
\def\dir{Dirichlet }
\def\Fa{{\oa F}}
\def\TR{{\cal T}}
 \def\ISO{{\rm ISO_p}}
 \def\Span{{\rm Span}}

\def\ALL{1}
\def\AV{2}
\def\BTA{3}
\def\BL{4}
\def\BRE{5}
\def\CNS{6}
\def\CP{7}
\def\CPW{8}
\def\CIL{9}
\def\CRA{10}
\def\DTT{11}
\def\DON{12}
\def\CG{13}
\def\DDD{14}
\def\DDR{15}
\def\GEO{16}
\def\HYP{17}
\def\BEL{18}
\def\SURVEY{19}
\def\AC{20}
\def\NOTES{21}
\def\AET{22}
\def\TANG{23}
\def\TANGG{24}
\def\LAG{25}
\def\SLE{26}
\def\JTY{27}
\def\KRY{28}
\def\PLI{29}
\def\RT{30}
\def\SLO{31}
\def\SPR{32}
\def\TRUU{33}
\def\TRU{34}
\def\TWC{35}
\def\TWCC{36}
\def\TWCCC{37}
\def\WAL{38}

\def\II{1}
\def\AA{2}
\def\AB{3}
\def\EE{4}
\def\BB{5}
\def\CC{6}
\def\CCC{7}
\def\DD{8}

\vskip .4in

\def\E{E}
\def\fpsi{{F_f(\psi)}}
\def\bL{{\bf \Lambda}}
\def\bdf{{\bf f}}
\def\UU{U}
\def\bbm{{\bf M}}
\def\gg{{\mathfrak g}}
\def\gra{\delta}
\def\Gm{G}
\def\LF{f}
\def\oN{^{1\over N}}

\font\headfont=cmr10 at 14 pt
\font\aufont=cmr10 at 11 pt


\title[DETERMINANT MAJORIZATION AND THE WORK OF GUO, PHONG AND TONG]
{\headfont  DETERMINANT MAJORIZATION AND THE WORK OF GUO-PHONG-TONG AND ABJA-OLIVE}

\date{\today}
\author{ F. Reese Harvey and H. Blaine Lawson, Jr.}
\thanks
{Partially supported by the Simons Foundation}

\maketitle

\vskip .5in

\centerline{\bf Abstract}

\medskip
The objective of this note is to establish the Determinant Majorization Formula
$F(A)^{1\over N} \geq \det(A)^{1\over n}$ for all operators $F$ determined by
an invariant G\aa rding-Dirichlet polynomial of degree $N$ on symmetric $n \times n$
matrices. Here {\sl invariant} means under the group O$(n)$, U$(n)$ or Sp$(n)\cdot {\rm Sp}(1)$
when the matrices are real symmetric, Hermitian symmetric, or quaternionic Hermitian symmetric
respectively. We also establish this formula for the Lagrangian Monge Amp\`ere Operator.
This greatly expands the applicability of the recent work of Guo-Phong-Tong 
and Guo-Pong for differential equations on complex manifolds.  It also relates to the work
of Abja-Olive on interior regularity.

Further applications to diagonal operators and to operators depending on the 
ordered eigenvalues are given.  Examples showing the preciseness 
of the results are presented. 
For the application to Abja-Olive's work, and other comments in the paper, we establish 
some results for G\aa rding-Dirichlet operators in appendices.  
One is an exhaustion lemma  for the G\aa rding cone.
Another gives  bounds for higher order derivatives,  which result from their elegant expressions as functions of the 
G\aa rding eigenvalues.
 There is also  a discussion of the crucial assumption of the Central Ray Hypothesis.

 \vfill \eject


\ 
\vskip1in
\centerline{\bf Table of Contents} \bigskip


 \hskip .5 in  1. Introduction

 \hskip .5 in  2. The Basic Lemma -- Majorization of the Determinant

 \hskip .5 in  3. G\aa rding-Dirichlet Polynomial Operators

 \hskip .5 in  4. Invariant G-D Operators

 \hskip .5 in  5. Relevance to the Work of Guo-Phong-Tong and Guo-Phong

 \hskip .5 in  6.  Interior Regularity -- Work of Abja-Olive

 \hskip .5 in  7. On the Enormous Universe of Invariant G-D Operators

 \hskip .5 in  8. Determinant Majorization for the Lagrangian Monge-Amp\`ere 
 
  \qquad\qquad \ Operator

 \hskip .5 in  9.  Being a G\aa rding-Dirichlet Operator is not Enough for 
 
  \qquad\qquad \    Determinant  Majorization

 \hskip .5 in  10. Other Applications of the Basic Lemma

 \hskip 1in Appendix A.  The Interior of the Polar Cone
 
 \hskip 1in   Appendix B. The G\aa rding  Gradient Map
  
\hskip 1in  Appendix C.  The Exhaustion Lemma

 \hskip 1in Appendix D.  The G\"uler Derivative Estimates
 
 \hskip 1in   Appendix E. The Central Ray Hypothesis

\vfill\eject

\noindent
{\headfont  1. Introduction.}  In a recent fundamental paper [GPT], B. Guo, D. H. Phong, and F. Tong
established an {\sl apriori} $L^\infty$ estimate for the complex Monge-Amp\`ere equation on   K\"ahler manifolds
by purely PDE techniques.  Their main theorem applied to many other equations as well.
The basic hypothesis in the  theorem can be established by proving a determinantal majorization inequality
for the operator.  In a later paper Guo and Phong proved  many further results  
on Hermitian manifolds with the same hypothesis [GP$_1$] (see also  [GP$_2$]).
The point of this paper is to {\sl prove the determinantal majorization inequality
for every invariant G\aa rding-Dirichlet operator,} which means that all the results above apply to this very large family of differential equations.  This inequality also appears as a hypothesis in work of S. Abja, S. Dinew and G. Olive [ADO],  [AO] on regularity, where there is also a compactness hypothesis which we establish in some generality.  So their results also apply to this same large constellation of equations.

The operators in question can be defined as follows.  Let $p(\l_1, ... , \l_n)$ be a real homogeneous polynomial
of degree $N$ on $\rn$ which is symmetric in the variables $\l_1, ... , \l_n$.  We assume that all the coefficients of $p$
are $\geq 0$.
 
We now define a polynomial operator $F:\Symn \to \bbr$,  of degree $N$ on symmetric $n\times n$-matrices, by
$$
F(A) \ \equiv \ p(\l_1(A), ... , \l_n(A))
\eqno{(1.1)}
$$
where $\l_1(A), ... , \l_n(A)$ are the eigenvalues of $A\in\Symn$.
Note that $p(\l_1(A), ... , \l_n(A))$ only makes sense because the  polynomial $p(\l_1, ... \l_n)$
 is {\bf symmetric}, i.e.,    invariant under permutations of the $\l_j$'s. 
 Since  $\l(A) = \{\l_1(A), ... , \l_n(A)\}$ is invariant under  
  the action of O$(n)$ acting by conjugation on $\Symn$, $F(A)$ is also O$(n)$-invariant.
This means that $F$ {\sl defines an operator on every riemannian manifold} by using the riemannian hessian [HL$_2$].

We can do the analogous thing in the complex case.  Let $A_\bbc$ be a hermitian symmetric $n\times n$-matrix
with eigenvalues $\l_k(A_\bbc)$ for $k=1, ... , n$.  Taking $p$ as in the above real case, we can define
$$
F(A_\bbc) \ \equiv \  p(\l_1(A_\bbc), ... , \l_n(A_\bbc)).
\eqno{(1.2)}
$$
This operator in invariant under the unitary group U$(n)$ and {\sl makes sense on any Hermitian complex
(or almost complex) manifold} by using $i\partial\dbar u$ and the metric.

Similar remarks can be make in the quaternionic case.

In these cases we shall prove the determinantal majorization inequalities:
$$
\begin{aligned}
F^{1\over N}(A) \ &\geq \ F(I)^{1\over N} \, \det(A)^{1\over n}  \\
F^{1\over N}(A_\bbc) \ &\geq \ F(I)^{1\over N} \, \det_\bbc(A_\bbc)^{1\over n}  \\
F^{1\over N}(A_\bbh) \ &\geq \ F(I)^{1\over N} \, \det_\bbh(A_\bbh)^{1\over n}
\end{aligned}
\eqno{(1.3)}
$$
where $F(I)>0$.  All three follow immediately from Basic Lemma 2.1.

We shall also prove this inequality for the Lagrangian Monge-Amp\`ere operator ${\rm M}_{\rm Lag}$
which is also defined on any almost complex  Hermitian manifold (in particular, on sympletic manifolds with a Gromov metric).
 The defining polynomial is unitarily invariant,
but the proof here is different because  this polynomial is zero on all traceless Hermitian matrices,
and so it is out of the category of operators considered in (1.3).
In Section 8 we prove that in complex dimension $n$,
$$
{\rm M}_{\rm Lag}(A)^{1\over 2^n} \ \geq \ \det_{\bbr}(A)^{1\over 2n}.
\eqno{(1.4)}
$$
See [HL$_5$] for complete details concerning this new operator.

The inequalities (1.3) also hold for polynomials $p$ which are not symmetric but satisfy the Central Ray Hypothesis in Section 2.
One just  orders the eigenvalues $\l_1(A)\leq \cdots \leq \l_n(A)$ and defines $F(A)$ by (1.1).  Here $F$ is not 
necessarily a polynomial.  (See Section 10.)

Determinant majorization is also established for {\sl diagonal operators}  
$$
p(a_{11}, a_{22}, ... , a_{nn})
$$
 where $p$ satisfies 
the hypotheses of Basic Lemma 2.1. 
However, if one drops the Central Ray Hypothesis, determinant majorization can fail.
We show that for examples of this type local interior regularity can also fail.  
This is done in  Section 9 .

As mentioned before the family of polynomial operators $F$ for which the determinant majorization holds is huge.
One way to see this is to consider  G\aa rding-Dirichlet (or G-D) operators on $\Symn$, which are invariant. 
These are homogeneous real polynomials $F$  on symmetric matrices with $F(I)>0$  such that: 

(i) $t\mapsto F(tI+A)$ has all real roots for every $A$, 

(ii) The {\sl G\aa rding cone} $\G$ contains $\{A : A>0\}$, where $\G$ is defined as the connected component of $\Symn-\{F=0\}$ containing 
the identity $I$,

(iii) $F$ is invariant under the action of O$(n)$ by conjugation on $\Symn$.

This family contains all elementary symmetric functions,  the $p$-fold sum operator (Example 3.3), and many more (see [HL$_3$, \S 5]).  
The set of these operators is closed under products, under  directional derivatives  in G\aa rding cone directions,
 and under a certain "composition" rule (see Section 7).

Every such operator $F$ can be expressed, in terms of the eigenvalues $\l(A)$ of $A$, as $F(A) = p(\l(A))$
where $p(\l_1, ... , \l_n)$ is a 
symmetric homogeneous polynomial which is $>0$ on $\rn_{>0}$. After nomalizing  by a positive constant, $p$ satisfies the 
three  Assumptions of the Basic Lemma 2.1. 
Assumption (1) follows from this 
positivity, and Assumption (3) follows from the O$(n)$-invariance. 
This  polynomial $p$  transfers over directly to the complex and quaternionic cases by
applying it to the eigenvalues of $A_\bbc$ or $A_\bbh$, and {\bf so (1.3) holds.}

We want to point out that   all the operators considered in this paper are elliptic.
In particular if $F$ is a G-D polynomial operator with G\aa rding cone $\G$, then 
a strengthened form of weak  ellipticity, namely
$\bra {\nabla F} P > 0$ for all $P> 0$, holds in its domain $\G$.

For the applications to the work of [ADO or AO] we need an Exhaustion Lemma which is proved in Appendix C.
For the log of  G-D polynomials there are elegant higher order
estimates which we thought would be good to include here (Appendix D).   An important assumption in the Basic Lemma
is the Central Ray Hypothesis, which is discussed from several points of view in Appendix E.
In a forthcoming article we shall give  a comprehensive introduction to  G-D operators and G\aa rding theory.
However, in this paper the focus in centered on the following  Basic Lemma 2.1.

We want to express our gratitude to Guillaume Olive for his careful reading of and helpful remarks on
the first version of this paper.

\vfill \eject

\noindent
{\headfont  2. The Basic Lemma -- Majorization of the Determinant}

Let $p$ ($\not \equiv 0$) be a real homogeneous polynomial of degree $N$ on $\rn$, and set $e\equiv (1,1,...,1)\in \rn$.
We write
$$
p(x) \ =\  \sum_\a a_{\a_1 \cdots \a_n} x_1^{\a_1}\cdots x_n^{\a_n} \ =\ \sum_{|\a|=N} a_\a x^\a.
$$

\noindent
{\bf Assumptions:}

\noindent
(1) \ \ All coefficients $a_\a \geq 0$.

\noindent
(2)   \ \ (Normalization)

$$
p(e) \ =\ \sum_\a a_\a \ =\ 1.
$$

\noindent
(3) \ \ The crucial Central Ray Hypothesis:
$$
D_e p \ =\ ke \qquad {\rm for}\ \ k>0,
$$
i.e., 
$$
{\partial p\over \partial x_1}(e) \ =\ \cdots\ = {\partial p\over \partial x_n}(e) \ =\ k >0.
$$ 
\noindent
Note that 
$$
{\partial p\over \partial x_j} (x) =\ \sum_\a a_\a \a_j x_1^{\a_1} \cdots x_j^{\a_j-1} \cdots x_n^{\a_n} 
\eqno{(2.1)}
$$
so that
$$
{\partial p \over \partial x_j} (e) \ =\ \sum_\a a_\a\,\a_j,  \qquad \text{ for each $j=1,...,n$}.
\eqno{(2.2)}
$$
Therefore, (3) can be restated as 

\noindent
(3)$'$   \ \  The  Central Ray Hypothesis:
$$
\sum_{|\a|=N} a_\a \a_j \ =\ k \qquad\text{for each $j=1, ... , n$ with $k>0$.} 
$$

\vskip .15in

It turns out that the following lemma goes back to  
a paper of Leonid Gurvits [Gur] in combinatorics.

\vskip.2in

\noindent
{\bf BASIC LEMMA 2.1.}  {\sl   For all $x_1>0, ... , x_n>0$,}
$$
p(x)^{1\over N}   \ \geq \ (x_1\cdots x_n)^{1\over n}.   
$$


\pf We first note that:
$$
\begin{aligned}
{1\over N} \log \, p(x) \ &=\  {1\over N}  \log \left\{   \sum_{|\a|=N} a_\a x^\a \right\} \\
&\geq  \  {1\over N}  \sum_{|\a|=N} a_\a \log (x^\a)   \quad\text{by concavity of log,  and (1) and (2)} \\
&= \  {1\over N}  \sum_{|\a|=N}  a_\a \{\a_1 \,\log x_1 +\cdots +\a_n\log\, x_n\}  \\
&= \ {1\over N} \sum_{j=1}^n \left\{  \sum_{|\a|=N}  a_\a \a_j \log\, x_j\right\} 
= \ {k\over N} \sum_{j=1}^n \log\, x_j \qquad \text{by (3)$'$}  \\
&= \ {k\over N} \log\, (x_1\cdots x_n).
\end{aligned}
$$

The value of $k$ is determined by Euler's Formula
$$
\bra {D_x p}  x \ =\ N p(x),
$$
followed by  assumptions (3) and then (2), yielding
$$
\bra {D_e p}  e\ =\ \bra {ke} e\ =\ kn\ =\ N p(e) \ =\ N,
$$
i.e.,
$$
k\ =\ {N\over n}.
$$
Substituting $1/n$ for $k/N$ in the inequality above, we have
$$
{1\over N} \log\, p(x) \ \geq\  {1\over n}  \log\,(x_1\cdots x_n).  \hskip .3in \mathqed
$$

\noindent
{\bf Corollary 2.2.} {\sl
Let $p(\l)= p(\l_1, ... , \l_n) \not \equiv 0$ be a symmetric homogeneous (real) polynomial of degree $N$
all of whose coefficients are $\geq 0$. Normalize $p$ so that $p(1,...,1)=1$. Then}
$$
p(\l_1, ... , \l_n)^{1\over N} \ \geq\  (\l_1\cdots \l_n)^{1\over n}
\eqno{(2.3)}
$$
for all $\l \in \rn_{>0}$.

\pf
We need to prove the Central Ray Hypothesis, that all the ${\partial p\over \partial \l_j}(e)$ agree (and equal $k$). This  follows from the invariance under
permutations.  Then, as noted  above, $k=N/n >0$.\qed

This corollary can be restated on the space $\Symn$ of second derivatives.


\noindent
{\bf Theorem 2.3.} {\sl
Let $F$ be a homogeneous polynomial of degree $N$ on $\Symn$ given by
$$
F(A) \ \equiv\ p(\l_1, ... , \l_n)
\eqno{(2.4)}
$$
where $p$ is as in the Corollary 2.2 and $\l_1, ... , \l_n$ are the eigenvalues of $A$.
Note that $F$ is normalized so that $F(I)=1$.
Then 
$$
F(A)^{1\over N}   \ \geq \  (\det\, A)^{1\over n}.
\eqno{(2.5)}
$$
for all $A\geq 0$.
}

\noindent{\bf Note 2.4.}
If we drop the normalization Assumption (2) on $p$,  and then replace $p(x)$ by $p(x)/p(e)$ which satisfies (2), we have for all $A\geq 0$ that
$$
  F(A)^{1\over N}  \ \geq\  \g\,  (\det\, A)^{1\over n}, \quad{\rm where} \ \g \ \equiv \  F(I)^{1\over N}\ =\ p(e)^{1\over N}.
\eqno{(2.6)}
$$

We will apply Theorem 2.3 to real symmetric, complex Hermitian symmetric,  quaternionic Hermitian symmetric matrices $A$,
and also the $2\times 2$ octonian Hermitian symmetric matrices.

\vskip.3in

\noindent
{\headfont 3.  G\aa rding-Dirichlet Polynomial Operators.} 

 Let $F:\Symn \to \bbr$ be a real  homogeneous polynomial of degree $N$
which is {\sl G\aa rding hyperbolic with respect to the Identity} $I$.  This mean that for each $A\in \Symn$, the polynomial
$t\mapsto F(tI+A)$ has all real roots.  The {\sl G\aa rding cone}  $\G$ of $F$ is the connected component of $\Symn-\{F=0\}$
which contains $I$. This is an open convex cone with 0 as the vertex. and for all directions $B\in\G$, $F$ is $B$-{\sl hyperbolic},  which
means that the polynomial $t\mapsto F(tB+A)$ has all real roots for every $A\in \Symn$. 
We always assume that $F(I)>0$. For complete details, see [G\aa r], [HL$_{3,4}$].

\noindent
{\bf Definition 3.1.}   We assume a {\bf positivity hypothesis}: $\G\supset \{A: A > 0\}$. This defines
a {\bf  G\aa rding-Dirichlet polynomial}, and applying this to the second derivative $D^2 u$  of functions $u$  in $\rn$ gives a {\bf  G\aa rding-Dirichlet (or G-D) operator.}  We shall always discuss a G-D operator in terms of the G-D polynomial that defines it.  So the terms "operator" and "polynomial" will be used interchangeably here.

These operators  have a form of strict ellipticity  on the G\aa rding cone $\G$:
$$
F(A+P) > F(A) \quad \text{for all  $A\in \G$ and $P>0$.}
\eqno{(3.1)}
$$
Recall that $F^{1\over N}$ is concave in $\G$ (see [G\aa r] or (D.2) below), which implies
$$
F(A+P)^{1\over N}  \geq F(A)^{1\over N} + F(P)^{1\over N} \quad  \text{for all  $A\in \G$ and $P>0$.}
\eqno{(3.2)}
$$
Combining this with (2.5) gives a stronger form of (3.1):
$$
F(A+P)^{1\over N} - F(A)^{1\over N}   \ \geq \  (\det\, P)^{1\over n}\quad\ \text{for all  $A\in \G$ and $P>0$.}
\eqno{(3.3)}
$$
where $F(I)=1$.  Letting $A\to 0$ in (3.3) yields (2.5), so (3.3) and (2.5) are equivalent  [ADO].

Basic G-D operators include the elementary symmetric functions $\s_k(A)$ (or Hessian operators).
Two new ones are the $p$-fold sum operator (3.5) below, and the Lagrangian Monge-Amp\`ere operator mentioned above
 (see [HL$_3$] for these and many others including for example
the  $\d$-uniformly elliptic operators $Q_\d(A) = \prod_j (\l_j(A) + \d \tr \l(A))$).

We also have the following construction.

\noindent
{\bf Proposition 3.2.}  {\sl  Let $F$ be a G-D operator of degree $N$.  Then for any $P>0$ the operator 
$$
\gra_P F(A) \ \equiv\ {d\over dt} F(tP+A) \biggr|_{t=0} \ = \ \bra {(\nabla F)_A} P
$$
is a G-D operator of degree $N-1$.
}

\pf The zeros of $\vf'(t)$ interlace the $N$ real zeros of $\vf(t) \equiv F(tP+A)$, so
$\gra_P F$ is G\aa rding. 
This also implies that the G\aa rding cones satisfy   $\G_{\gra_P F} \supset \G_F$ (see [G\aa r]), and 
$\G_F \supset \{A: A>0\}$, so $\gra_P F$ satisfies positivity.
  \qed

This construction can be iterated.
  For a G-D operator $F$, let $F(tI+A) = \prod_{j=1}^N (t+\l_j^F(A))$ 
  where $\l_j^F(A)$ are called the {\bf eigenvalues} of $F$ with respect to $I$.
  Then differentiating this product at $t=0$ and applying Proposition 3.2 yields that
  $$
  {1\over k!} \gra^k_I F(A) \ =\  \s_{N-k}(\l_1^F(A), ... , \l_n^F(A))
  \eqno{(3.4)}
$$
for $k=0,1, ... , N$ are G-D polynomials, where  $\s_j$ is the $j^{\rm th}$ elementary symmetric function.
It follows that {\sl the elementary symmetric functions of $\l_1^F(A), ... , \l_n^F(A)$ are polynomials in $A$.}

The next example was introduced  in [HL$_1$], along with $p$-geometry/$p$-potential theory.

\noindent
{\bf Example 3.3. (The $p$-Fold Sum Operator).}  Let 
$$
F(A) = \prod_{|J|=p} \l_J(A)
\eqno{(3.5)}
$$
 where 
$$
\l_J(A) \ =\ \l_{j_1}(A) + \cdots + \l_{j_p}(A)\qquad{\rm where}\ \  j_1< \cdots < j_p \ \ {\rm and\ degree =}\ \ N={n\choose p}.
$$
(The inequality (2.5) for this example was first proven in [D].)

We shall see in \S 7 that the family of invariant G-D polynomials is huge! Concerning ellipticity.

\noindent
{\bf Proposition 3.4.}  {\sl A complete G-D operator $F$ has uniformly  elliptic linearization at each point of its G\aa rding cone $\G$.
That is, at $A\in\G$
$$
B \mapsto \text{$\bra {(\nabla  F)_A}{B}$ has coefficient matrix $(\nabla F)_A > 0$.}
$$
}
\pf  This is  Proposition B.1 in Appendix $B$.
\qed

Note that the  uniform ellipticity $(\nabla F)_A>0$ is stronger than 
$\bra {(\nabla F)_A}{P} >0$, $\forall \,P>0$ mentioned earlier, and it requires the completeness hypothesis.
Completeness means that all the variables in $\rn$ are needed to define $F$.
This rules out operators such as $F(D^2 u) = \partial^2 u/\partial x_1^2 \ (n\geq 2)$.
The invariant operators in the next section are always complete.


\vskip.3in

\noindent
{\headfont  4. Invariant G-D Operators.}

We shall be concerned principally with three cases.

\noindent
{\bf Case 1.  (Real Invariant G-D Operators.)}  A G-D polynomial $F$ on $\Symn$ of degree $N$ is 
{\bf real invariant} if 
$$
F(A) \ = \ F(gAg^t) \qquad \text{for all $g\in {\rm O}(n)$}.
$$
Every such polynomial canonically determines an operator on every riemannian manifold.

By Theorem 4.1  below, $F$ {\bf  is real invariant if and only if }
$$
F(A) \ =\ p(\l_1(A), ... , \l_n(A))
$$
where $p$ is a homogeneous symmetric polynomial of degree $N$ in the eigenvalues of $A$.
This polynomial $p(\l_1, ... , \l_n)$ on $\rn$ satisfies the Assumptions (1) and (3) of the Basic Lemma 2.1,
and $p(e) = F(I) >0$.

\noindent
{\bf Case 2.  (Complex Invariant G-D Operators.)} 
 Let $\Sym_\bbc(\bbc^n)$ be the space of Hermitian symmetric $n\times n$-matrices on $\bbc^n$.
A (real-valued) G-D polynomial $F$ on $\Sym_\bbc(\bbc^n)$ of degree $N$ is  {\bf complex invariant} if 
$$
F(A_\bbc) \ = \ F(gA_\bbc \overline g^t) \qquad \text{for all $g\in {\rm U}(n)$}.
$$

 Every such polynomial $F$ canonically determines an operator on every Hermitian complex (or almost complex) manifold.

By Theorem 4.1 below, $F$ {\bf  is complex invariant if and only if }
$$
F(A_\bbc) \ =\ p(\l_1(A_\bbc), ... , \l_n(A_\bbc))
$$
where $p$ is a homogeneous symmetric polynomial of degree $N$ in the eigenvalues of $A_\bbc$.
This polynomial $p(\l_1, ... , \l_n)$ on $\rn$ satisfies the Assumptions (1) and (3) of the Basic Lemma 2.1,
and $p(e) = F(I) >0$.

This case can be looked at from the real point of view.
 Let $\bbc^n = (\bbr^{2n}, J)$.  Then $\Sym_\bbc(\bbc^n)$ can be identified with the subset of $A\in\Sym(\bbr^{2n})$ satisfying $AJ=JA$.
  There is a projection $(\cdot)_\bbc:\Sym(\bbr^{2n}) \to \Sym_\bbc(\bbc^n)$ given by setting 
$$
A_\bbc \ =\ {1\over 2} (A-JAJ)
$$
The real eigenvalues of $A_{\bbc}$ are the eigenvalues of $A_{\bbc}$ viewed as a Hermitian symmetric matrix,
but where each eigenvalue now appears {\bf twice}.

\noindent
{\bf Case 3.  (Quaternionic Invariant G-D Operators.)} 
Here we take the real point of view.
Let $\bbh^n = (\bbr^{4n}, I, J, K)$ where $I,J,K$ are 
 complex structures satisfying the usual relations.
 (Thus $\bbh^n$ can be viewed as a right  quaternion vector space with scalars $\a+\b I + \g J + \d K$.)
 Then $\Sym_\bbh(\bbh^n)$, the space of quaternionic Hermitian  symmetric $n\times n$-matrices, can be considered as the subset of
$A\in \Sym(\bbr^{4n})$ which commute with $I,J$ and $K$.  

  There is a projection $(\cdot)_\bbh:\Sym(\bbr^{4n}) \to \Sym_\bbh(\bbh^n)$ given by setting 
$$
A_\bbh \ =\ {1\over 4} (A- IAI - JAJ - KAK).
$$
Using the natural inner product on $\Sym(\bbr^{4n})$, this map from $A$ to $A_\bbh$ is orthogonal projection.  
Considering $A_\bbh \in \Sym_\bbh(\bbh^{n})$ as an element of $ \Sym(\bbr^{4n})$,  it has a canonical form under the conjugate
action of  O$(4n)$, namely
$$
A_\bbh \ =\ \l_1 P_{W_1} + \cdots + \l_s P_{W_s}
\eqno{(4.1)}
$$
with $\l_1< \l_2 < \cdots < \l_s$ and $\bbr^{4n} = W_1 \oplus \cdots \oplus W_s$ an orthogonal decomposition, comprising
the distinct eigenvlaues and corresponding eigenspaces.  Here $P_W$ denotes orthogonal projection onto the real subspace
$W\ss\bbr^{4n}$.  Moreover, if $v$ is an eigenvector of $A_\bbh  \in \Sym(\bbr^{4n})$ with eigenvalue $\l$, then $vI, vJ, vK$ are also 
eigenvectors with the same eigenvalue $\l$.  That is,
$$
\text{each eigenspace $W_j$ is a right quaternion vector subspace of $\bbh^n$}.
\eqno{(4.2)}
$$
Taken together, (4.1) and (4.2) provide the {\bf canonical form} for $A_\bbh \in \Sym_\bbh(\bbh^n)$
under the conjugate action of Sp$(n)\cdot {\rm Sp}(1) \equiv {\rm Sp}(n)\times_{\bbz_2} {\rm Sp}(1)$.
(Note that this canonical form could also be expressed using the eigenvalues $\l^0_j(A_\bbh), j=1, ... , n$, listed to multiplicity.)

A homogeneous real-valued polynomial $F$ on $\Sym_\bbh(\bbh^n)$ is G-D if $t\mapsto F(tI + A_\bbh)$ has all real roots
for every $A_\bbh \in \Sym_\bbh(\bbh^n)$.

Let $F$ be a  real homogeneous polynomial $F$ on $\Sym_\bbh(\bbh^n)$ of degree $N$.  Then $F$ is 
{\bf quaternionic invariant} if 
$$
F(A_\bbh) \ = \ F(gA_\bbh g^t) \qquad \text{for all $g\in {\rm Sp}(n) \cdot {\rm Sp}(1)$}.
$$
Every such polynomial canonically determines an operator on every manifold with a topological $\bbh$-structure and a compatible metric,
in particular a hyperK\"ahler manifold with a Calabi-Yau metric.

By Theorem 4.1 below, a G-D polynomial $F$ of degree $N$ on   $\Sym_\bbh(\bbh^n)$
 {\bf  is quaternionic invariant if and only if }
$$
F(A_\bbh) \ =\ p(\l_1^0(A_\bbh), ... , \l_n^0(A_\bbh))
\eqno{(4.3)}
$$
where $p$ is a homogeneous symmetric polynomial of degree $N$, 
and $\l^0_j(A_\bbh)$ are the eigenvalues of $A_\bbh$ listed to multiplicity.
This polynomial $p(\l_1, ... , \l_n)$ on $\rn$ satisfies the Assumptions (1) and (3) of the Basic Lemma 2.1,
and $p(e) = F(I) >0$.


\noindent
{\bf Theorem 4.1.} {\sl Let $F(A)$  be a real invariant G-D polynomial of degree $N$ on $\Symn$.
Then 
$$
F(A) \ =\  p(\l_1(A), ... , \l_n(A))
$$
where $p(\l_1, ... , \l_n)$ is a symmetric homogeneous polynomial of degree $N$ on $\rn$ which satisfies
Assumptions (1) and (3) of Lemma 2.1 and for which $p(e)>0$.  Furthermore, $p$ is G-D, that is 
$p$ is G\aa rding hyperbolic w.r.t.\ $e$ and the set $\{\l\in \rn:  \l_j>0\ \forall\, j\}$ is contained in its G\aa rding cone $\G(p)$.

Conversely,  if $F$ is given this
way for such a $p$, then $F$ is a real invariant G-D operator.

The analogous statements hold in the complex and quaternionic cases..
}

Before the proof we make several remarks.

\noindent
{\bf Remark 4.2.}
It is important to note  in the three determinant cases:
$$
(1)\ \ \det_\bbr A, \ A\in \Sym_\bbr(\bbr^n), \qquad 
(2)\ \ \det_\bbc A_\bbc, \ A\in \Sym_\bbr(\bbc^n), $$ $$ \quad {\rm or} \quad
(3)\ \ \det_\bbh A_\bbh, \ A\in \Sym_\bbr(\bbh^n), \qquad 
$$
that $F$ is a G-D polynomial with G\aa rding $I$-eigenvalues given by 
$$
(1) \ \ \l(A) \ =\ \{\l_1(A), ... ,, \l_n(A)\},  \qquad
(2) \ \ \l(A) \ =\ \{\l_1(A_\bbc), ... ,, \l_n(A_\bbc)\}, $$ $$\quad {\rm or}\ 
(3) \ \ \l(A) \ =\ \{\l_1(A_\bbh), ... ,, \l_n(A_\bbh)\}, {\rm respectively},
$$
defined in each case by the appropriate canonical form.

\noindent
{\bf Remark 4.3.}
The fact that $F(A) \equiv \det_\bbh A_\bbh \equiv \l_1(A_\bbh)\cdots \l_n(A_\bbh)$ is a polynomial in 
$A\in \Sym_\bbr(\bbh^n)$ (or equivalently a polynomial in $A_\bbh \in \Sym_\bbh(\bbh^n)\ss \Sym_\bbr(\bbr^{4n})$) requires proof.
The Moore determinant provides an algebraic construction which is rather complicated (see [AV].)
Some might prefer the following proof.

\noindent
{\bf Lemma 4.4.}  {\sl  The determinant 
$
\det_\bbh A_\bbh  \equiv \l_1(A_\bbh)\cdots \l_n(A_\bbh)
$
is a polynomial in $A \in \Sym_\bbr(\bbh^n)$.}

\pf
Let $\mathbb H = (\bbr^{4n}, I,J,K)$ and let $\pi:\Sym(\bbr^{4n}) \to \Sym(\bbr^{4n})$ be given by 
$$
\pi(A) \ =\  {1\over 4} (A - IAI-JAJ-KAK)\ \equiv\ A_\bbh.
$$
Consider 
$$
\begin{aligned}
\det(I+tA_\bbh) \ &=\ 1+\Sigma_1(A)t + \Sigma_2(A)t^2 + \cdots + \Sigma_{4n}t^{4n}   \\
&= \ (1+\s_1(\l)t + \s_2(\l)t^2 + \s_3(\l)t^3 + \cdots + \s_n(\l)t^n)^4
\end{aligned}
$$
where the $\l$'s are the $n$ eigenvalues of $A_\bbh$.
Now $\det(I+tA_\bbh)$ is a polynomial in $(A,t)$ so the $\Sigma_k(A)$'s are polynomials in $A$.
We also have
$$
\begin{aligned}
 &(1+\s_1(\l)t + \s_2(\l)t^2 + \s_3(\l)t^3 + \cdots) (1+\s_1(\l)t + \s_2(\l)t^2 + \s_3(\l)t^3 + \cdots) \cdot \\
 &\cdot  (1+\s_1(\l)t + \s_2(\l)t^2 + \s_3(\l)t^3 + \cdots) (1+\s_1(\l)t + \s_2(\l)t^2 + \s_3(\l)t^3 + \cdots) \\
 &=\  1 + 4\s_1(\l)t + [4\s_2(\l) + 6\s_1^2(\l)]t^2 + [4\s_3(\l) + \cdots]t^3 + [4\s_4(\l)+ \cdots]t^4 + \cdots \\
 &= \ \ 1+\Sigma_1(A)t + \Sigma_2(A)t^2 +  \Sigma_3(A)t^3 =\Sigma_4(A)t^4 +\cdots 
  \end{aligned}
$$
Now 

\noindent
$4\s_1(\l) = \Sigma_1(A)$ is a polynomial in $A$.  Therefore

\noindent 
 $4\s_2(\l) + 6\s_1^2(\l)= \Sigma_2(A)$ is a polynomial in $A$.  Hence $\s_2(\l)$ is a polynomial in $A$.
 
  \noindent
 In general
 $4\sigma_k(\l) + {\rm a \ polynomial \ in\ } \ A \ =\Sigma_k(A)$ and so $\s_k(\l)$ is a polynomial in $A$.\qed

 A polynomial $p(\l)$ which satisfies all the conditions in Theorem 4.1 is defined (Def. 7.3) to be a {\bf  universal G-D polynomial}.

 Theorem 4.1 is generalized in Theorem 7.2.
 There the G\aa rding $I$-eigenvalues used in Theorem 4.1 for one of the three determinants 
 $\det_\bbr, \det_\bbc, \det_\bbh$ are replaced by the G\aa rding $I$-eigenvlues of any G-D operator.

\noindent
{\bf Proof of Theorem 4.1}.
We begin by restricting $F$ to the subspace  $\bbd \ss \Symn$ of diagonal matrices.
For each $\l \in \rn$, let $M(\l) = \sum _i \l_i e_i\circ e_i$ denote the diagonal matrix with entries $\l_1, ... , \l_n$. 
Define
$$
p(\l) \ \equiv \   F(M(\l)) \qquad {\rm for}\ \ \l \in \rn.
$$
Then $p(\l)$ is a 
 homogenous polynomial of degree $N$.  The fact that $F$ is invariant 
under the orthogonal group implies that $p$ is invariant under permutations of the $\l_j$'s, i.e., $p$ is a symmetric
polynomial.  Now by the conjugation invariance, this polynomial determines $F$ by the formula $F(A) = p(\l_1(A), ... , \l_n(A))$.
So now everything is reconstructed in terms of $p$.  The positivity hypothesis that $F>0$ on $\{A:A>0\}$ is equivalent to the
strict positivity of $p$ on $\bbr^n_{>0}$, since $A>0 \ \iff\ \l(A) \in\rn_{>0}$.   Assertion (1) that all the coefficients of $p$ are $\geq 0$, can be proved as follows.

Note that
$$
{1\over \a_1! \cdots \a_n!}\, {\partial^{|\a|}  p \over \partial x_1^{\a_1}\cdots  \partial x_n^{\a_n}}
$$
equals the coefficient $a_\a$ for each multi-index $\a$ of length $|\a| \equiv \sum_k \a_k = N$.
To prove that this $\a^{\rm th}$  partial of $p$ is $\geq 0$ we use Proposition 3.2.

Since $p(\l) = F(\l_1 e_1\circ e_1 + \cdots + \l_n e_n\circ e_n)$, the partials of $p$ 
$$
{\partial p  \over \partial \l_j}(\l) \ =\ \bra {(\nabla F)_{M(\l)}}  {\ e_j \circ e_j}
$$
 are equal to the directional derivatives of $F$ in the directions $P_j = e_j\circ e_j \geq 0$
By Proposition 3.2   we have that  $\bra {(\nabla F)_A}  {P} = (\gra_P F)(A)>0$ for all $A\in \G$ and $P>0$.  Hence, 
$(\gra_P F)(A)\geq 0$ for all $A\in \G$ and $P\geq 0$.  

By Proposition 3.2, $\gra_P F$ is also an invariant G-D operator. Therefore, this process can be repeated $N$ times proving that 
the $N^{\rm th}$ directional derivatives in directions $P_1>0, ... , P_N>0$ are $>0$.  Hence, 
the $N^{\rm th}$ directional derivatives in directions $P_1\geq 0, ... , P_N\geq 0$ are $\geq0$.
Taking these latter $P_j$'s to each be in axis directions $e_k\circ e_k$ yields the desired result that the $\a^{\rm th}$ partial of $p$,
$|\a|=N$, is $\geq 0$.
So Assertion (1) is proved.

For Assumption (3) we use the fact that the gradient $D_e p$ at $e$ is invariant under permutation 
of the $\l_j$'s since $e=(1, ... , 1)$ is invariant and so
is $p$.   Since $\bbr e$ is the only subspace of $\rn$  fixed under the permutation group, we have $D_e p = ke$ for some $k\in \bbr$.
However $k$ must be positive since by Euler's formula
$
nk \ =\   \bra { D_ep} e  \ =\ Np(e) >0
$
as in the proof of Basic Lemma 2.1.

The arguments for the complex and quaternionic cases are exactly the same.
\qed

\noindent
{\bf Theorem 4.5.} {\sl
Let $F$ be an invariant G-D operator from one of the three cases above. Let $\det(A)$ denote the complex
or quaternionic determinant in the complex and quaternionic cases.  Then
$$
F(A)^{1\over N}  \ \geq \ F(I)^{1\over N} \det(A)^{1\over n}
\eqno{(4.4)}
$$
}
\pf
This follows directly from Theorem 4.1 and Basic  Lemma 2.1.\qed

\vskip .3in 


\noindent
{\headfont 5. Relevance to the  Work of Guo-Phong-Tong and Guo-Phong.}

One of the principle motivations for this work was the recent paper of 
B. Guo, D. H. Phong, and F. Tong [GPT], which among other things
gave a purely PDE proof of the $C^0$-estimate in Yau's proof of the Calabi Conjecture.
Their main theorem, Theorem 1, has much broader applicability to operators of complex invariant type
(as in Section 3 in this paper) on K\"ahler manifolds.  The main assumption (1.4) in this theorem follows  if one can prove the majorization
of the determinant formula in the cone $\Int \cp$ of positive definite matrices (see [GPT, Lemma 4]). That is exactly  what we have 
done for all complex invariant G\aa rding-Dirichlet operators, and, therefore, for all the associated operators on K\"ahler  manifolds.

In a later paper of Guo and Phong [GP], a number of important results for operators on Hermitian manifolds were established.
Again the important assumption is established by the majorization of the determinant formula. Hence,  all the
results in this paper hold for every operator  on a Hermitian manifold, which is induced by any complex invariant
G\aa rding-Dirichlet operator.

\bigskip


\noindent
{\headfont  6. Interior Regularity - Work of Abja-Olive.}   

 In this section we establish $C^\infty$ interior regularity for 
$W^{2,p}$-solutions of complex invariant G-D operators (in the Euclidean case).  This is based on results of Abja-Olive
[AO].  
We begin by showing that the basic assumptions (a) -- (f) in [AO] are satisfied.
As above we work in eigenvalue space with a  symmetric polynomial $p(\l_1, ... , \l_n)$ which satisfies the assumptions of Lemma 2.1. 
Let $\G$ denote the connected component of $\rn-\{p=0\}$ containing $e$.

(a)  \ \  $\bbr^n_{>0} \ \ss \ \G\ \ss\ \{\l_1+\cdots + \l_n >0\}$.

(b)  \ \ $p$ is a $C^1$ on $\bar \G$.

(c) \ \ $p>0$ on $\G$.

(d) \ \ $p$ is positively homogeneous of degree $N$ on $\G$.

(e) \ \ $p^{1\over N}$ is concave on $\G$.

(f) \ \ $p(\l)^{1\over N} \geq (\l_1\cdots \l_n)^{1\over n}$ with equality on span$(e)$.

\pf
The left hand inclusion of (a) follows from the  positivity hypothesis. For the right hand inclusion, note that 
$\G$ is convex and invariant under permutations of the $\l_j$'s. 
Among  half spaces $H$ containing $\G$ with
the vertex $0\in \partial H$, there is only  one $H_0$ which is  invariant under permutations. Since $\bbr\cdot e$ is the only
invariant line, we have $H_0 = e^\perp =\{\l_1+\cdots+\l_n\geq 0\}$.

Hypothesis $(b)$ is clear since $p$ is a polynomial.  Hypothesis $(c$) follows from the definition of $\G$.
Hypothesis (d) is an assumption.  Hypothesis (e) is Theorem 2 in [G\aa r] for the initial polynomial on matrices This restricts to concavity for diagonal matrices.  Hypothesis (f)  is  Corollary 2.2. 

Now the assumptions of [AO] are on matrix space and not on eigenvalue space $\rn$.  One needs to know that  $C\ss\rn$ is convex if and only if  $\l^{-1}(C)$ is convex. This was proved in [BGLS]. ( See also [HL$_4$, Thm. 8.4] and [AO, Prop. 3.1].) 
One also needs to know that if $C\ss\rn$ is compact, 
then so is $\l^{-1}(C)$. However, $C\times {\rm O}(n)$ is compact and $\l^{-1}(C)$ is an image of $C\times {\rm O}(n)$.
\qed

We can now apply the Main Theorem 1.1 in [AO] where everything takes place in  $\bbc^n$
and the authors deal with Hermitian symmetric matrices.
Let $F$ be the invariant G-D operator on $\bbc^n$  given by a homogeneous polynomial $p(\l_1, ... , \l_n)$ of degree $N$ as above,
applied to the eigenvalues of Hermitian symmetric matrices $A$, that is, $F(A)=p(\l(A))$.    We now let $\G$ denote the G\aa rding cone of $F$, and  we fix a domain $\O\ss\bbc^n$.

\noindent
{\bf Theorem 6.1.}  {\sl 
Let  $u \in W^{2,q}_{\rm loc}(\O)$, $q>n\, \max\{N-1,1\}$, with  $D^2u(x) \in\G$ for a.a.\ $x\in\O$ satisfying  
$$
F(D^2u) \ =\ f  \qquad {\rm a.e.\   in}\ \O
$$
 where $f >0$ and in $C^2(\O)$.
Then for each domain $\O_0\ss\ss\O$ there exists $R>0$ with 
$$
\|  \D u   \|_{ L^\infty(\O_0)} \ \leq \ R.
$$
}

We now apply this bound on $\D u$.  We begin with the following.

(g) \ \ (Sets of Uniform Ellipticity)  For each  $R>0$, $F$ is uniformly elliptic
in the region 
$$
\G_R \ =\ \left\{ A\in\G : \tr(A) < R \ \ {\rm and}\ \ {1\over R} < F(A)\right\}
$$

Now the polynomial $F$ has the property that for $A\in\G$,
$$
{D_A F} \in \G^* \ (\text{the "open polar" of $\G$}) \ss \Int \cp.
$$
The definition  of the   {\sl open polar}  $\G^*$ is in  Appendix A.  
One sees easily that an invariant D-G polynomial satisfies the completeness condition given in Appendix A.
Therefore, Proposition B.1 in Appendix B applies to show that ${D_A F} \in \G^*$.  Since $\Int\, \cp \ss \G$, we have $\G^*\ss \Int \, \cp$.
  Therefore the linearization of $F$ is positive definite on all of $\G$.   
 Assumption (g) is then a result of the following fact.
 
(g)$'$ \ \ For all $R>0$ the set
$$
 \overline{\G}_R \ \equiv\ \left\{ A \in \G : \sum_{j=1}^n \l_j(A) \leq   R \ \ {\rm and} \ \ {1\over R}  \leq F(A)  \right \}     \ \ \text{is compact}.   
$$

\pf  
From (a) above one easily sees  that $I\in \G^0$.
The assertion (g)$'$ then follows  from Corollary C.3  in Appendix C with $c_1 = R, c_2 = \log\,(1/R)$ and $\gg=F$.\qed

The fact that (g) holds in this generality is quite useful, so we state it separately.

\noindent
{\bf Proposition 6.2.}  {\sl
For each  $R>0$, the G-D operator $F$ restricted to the  set $\overline{\G}_R$ is uniformly elliptic.}

In particular, we have interior regularity for $W^{2,p}$ solutions (see [CC], [AO], [W]). 

\noindent
{\bf Theorem 6.3.} {\sl  Let $F$ be an invariant complex G-D operator.  Under the assumptions of Theorem 6.1, 
$$
\text{  $f \in C^\infty(\O) \ \ \Rightarrow \ \ 
u\in C^\infty(\O)$.}
$$
}

\vskip .3in

\noindent
{\headfont  7. On the Enormous Universe of Invariant  G-D Operators.} 

Let's fix one of the algebras $\bbr, \bbc$ or $\bbh$, and for the discussion in this 
chapter "invariant" will always mean "real invariant", "complex invariant" or "quaternionic invariant"
accordingly.  The discussion is identical in these three cases.

We begin with the following.

\noindent
{\bf Proposition 7.1.} {\sl
The product of two invariant G-D polynomials is again an invariant G-D polynomial, with eigenvalue set the union of the $F$ and $G$ eigenvalues.
}

\pf
If 
$$
F(tI+A) \ =\  F(I) \prod_{k=1}^n (t+ \l^F_k(A)) 
\ \ {\rm and}\ \ 
G(tI+A)\ =\ G(I) \prod_{k=1}^n (t+ \l^G_k(A)),
$$
then 
$$
F(tI+A)  G(tI+A) \ =\ F(I)G(I) \prod_{k=1}^n (t+ \l^F_k(A)) \prod_{k=1}^n (t+ \l^G_k(A)),
$$
and so $FG$ is G\aa rding hyperbolic w.r.t.\ $I$.  If $F$ and $G$  are each invariant under a group ${\mathcal G}$
 acting on the $A$'s, so is the product.\qed

\noindent
{\bf Proposition 7.2} 
{\sl    Let $F$ be an invariant G-D polynomial of degree $N$.  Then  the operator 
$$
(\gra_I F)(A) \ \equiv\ {d\over dt} F(tI+A) \biggr|_{t=0} \ =\ \bra{(\nabla F)_A}{I}
\eqno{(7.1)}
$$
is an invariant G-D polynomial of degree $N-1$.}

\pf  The invariance of $\gra_I F$ is clear.  The rest is Proposition 3.2.\qed

We formalize  the definition given after Theorem 4.1, which has been taken from [HL$_3$].

\noindent
{\bf Definition 7.3.}
A {\bf universal G-D polynomial} is a  homogeneous symmetric polynomial $p(\l_1,  ... ,\l_n)$
which is  G\aa rding hyperbolic w.r.t.\ $e$ and satisfies Assumptions (1) and (3) of Basic Lemma 2.1
with $p(e)>0$. Note that Assumption (1) implies $p>0$ on $\rn_{>0}$.

These are exactly the polynomials such that $F(A) = p(\l^F(A))$ for  invariant 
G-D polynomials $F$  (Theorem 4.1).

Consider now  a G-D polynomial  $F$ of degree $M$ on $\Symn$ with I-eigenvalues $\l^F(A) = (\l_1^F(A), ... , \l_M^F(A))$.
We shall construct a new  G-D polynomial  by  applying  a universal G\aa rding polynomial $p$ of degree N 
in $M$ variables to the  $F$-eigenvalues. 
This is why $p$ is called universal.

\noindent
{\bf Theorem 7.4.}  
{\sl
 Let $p$ and $F$ be as above.  Define 
$$
P_F(A) \ \equiv\ p(\l_1^F(A), ... , \l_M^F(A)) \ =\ p(\l^F(A))
$$
for $A\in \Symn$.  Then $P_F$ is also a G-D polynomial of degree $N$ on $\Symn$, with G\aa rding eigenvalues
$$
\l_j^{P_F}(A) \ =\ \L_j(\l^F(A))  \qquad j=1, ... , N,
\eqno{(7.2)}
$$  
where $\L_1(\l), ... , \L_N(\l)$, are the $e$-eigenvalues of the G-D polynomial $p$. \\
The G\aa rding cones are related by $\G_{P_F} = (\l^F)^{-1}(\G_P)$.

Furthermore, if $F$ is invariant, so is $P_F$.
}

\pf
We first observe that $P_F(A)$ is a polynomial.  This is evident if $p(\l) = \l_1\cdots \l_M$ where $P_F = F$.
 It then follows for all the elementary symmetric functions $\s_1, ... , \s_M$ by Proposition 3.2 and Example 3.3.
 Now by a classical result, $p(\l)$ is a polynomial in the elementary symmetric functions, and so 
 $P_F(A)$ is a polynomial.
 
Note first that $P_F(I) = p(\l^F(I)) = p(e) >0$.      By (2.1) in [G\aa r],
$$
p(te+\l^F) \ =\ p(e) \prod_{j=1}^N(t+ \L_j(\l^F)).
$$
Now $\l^F(I)=e$ since $F(tI+I) = (t+1)^MF(I)$. Therefore (2)$'$ in  [G\aa r] says that
$$
\l^F(tI+A) \ =\ te + \l^F(A),
$$
and so
$$
\begin{aligned}
P_F(tI+A) \ &=\ p(\l^F(tI+A))  \\ &=\ p(te+\l^F(A)) \\ &=\ p(e) \prod_{j=1}^N (t+\L_j(\l^F(A)).
\end{aligned}
$$
Therefore, $\L_j(\l^F(A))$, $j=1, ... , N$,  are the G\aa rding $I$-eigenvalues of $P_F$.
The statement about G\aa rding cones follows from (7.2). 
The invariance statement follows from Theorem 4.1. \qed

Many basic invariant G-D polynomials are given in [HL$_{3,4}$].  Taking these and 
using  Propositions 7.2, 7.3, and 7.4 successively in long chains show that the {\bf set of 
invariant G-D polynomials  operators is  very large.}

\def\e{e}

\vskip .3in

\noindent
{\headfont  8. Determinant Majorization for the Lagrangian Monge-Amp\`ere Operator.}  

The Lagrangian Monge-Amp\`ere Operator is a U$(n)$-invariant G-D operator on $\bbc^n$ which is dependent on 
the Skew part and independent of the traceless Hermitian symmetric part of a matrix $A$. Its
pointwise linearizations are elliptic. For a compete discussion, including its importance for a lagrangian
potential theory and its geometrical interpretation,   one should see [HL$_5$], where the operator was introduced.
Note that by its U$(n)$-invariance this operator  is defined on
any (almost) complex Hermitian manifold.  

For the algebra pertaining to this equation, we shall be  working on 
\newline $\bbc^n = (\bbr^{2n}, J)$. 
The operator is defined on the U$(n)$-invariant subspace 
$\bbr I \oplus {\rm Herm}^{\rm Skew}(\bbc^n) \ss\Sym_\bbr(\bbr^{2n})$ 
where $ {\rm Herm}^{\rm Skew}(\bbc^n) = \{A : JA=-AJ\}$. It is zero on the space
of traceless Hermitian symmetric matrices (those traceless $A$ with $JA=AJ$). Therefore we can  assume that 
this part of $A$ is zero.  If $Ae=\l e$ and $A$ is skew Hermitian, then $A(Je) = - \l Je$. Hence, we can 
let $\pm\l_1, ...  , \pm \l_n$ denote the eigenvalues of the skew part $A^{\rm skew} \equiv \half(A+JAJ)$. Then
$$
A \ =\ tI + A^{\rm skew}
$$
with eigenvlaues
$$
t\pm \l_1, ... , t\pm\l_n \qquad{\rm and} \qquad \tr(A) \ =\ 2nt.
$$
Thus we have
$$
\det(A) \ =\ \prod_{j=1}^n (t+\l_j)(t-\l_j).
$$
With $\mu \equiv {1\over 2} \tr(A) = nt$ the operator is defined in [HL$_5$]  by
$$
\begin{aligned}
F(A) \ &=\ \prod_{2^n \pm} (\mu\pm\l_1\pm \l_2\pm \cdots\pm \l_n) \ =\ \prod_{2^n  \pm} ( nt \pm\l_1\pm \l_2\pm \cdots\pm \l_n) \\
&=\ \prod_{2^n \pm}(\e_1^{\pm} + \cdots + \e_n^\pm)  \qquad {\rm where} \ \ \e_j^\pm \ =\ t\pm \l_j.
\end{aligned}
$$
$$
\text{The eigenvalues of $A$ are $ \e_j^+$ and $\e_j^-$ for $j=1,..., n$, so }
$$
 $$
\det(A)\ =\ \e_1^+ \e_1^- \cdots \e_n^+\e_n^-
$$

\noindent
{\bf Proposition 8.1}.  {\sl In the region where $A>0$, }
$$
F(A)^{1\over 2^n} \ \geq\ \det(A)^{1\over 2n}
$$
\pf
Positivity implies that  $\e_j^{\pm} >0$ for all $j$. Hence, for all $2^{n-1}$ choices of $\pm$ we have 
$$
\log(\e_1^\pm + \cdots +\e_j^+ + \cdots + \e_n^\pm) + \log(\e_1^\pm + \cdots +\e_j^- + \cdots + \e_n^\pm) 
\ \geq\ \log(\e_j^-) + \log(\e_j^-).
$$
Summing over all  $2^{n-1}$ choices gives
$$
\begin{aligned}
&{1\over 2^{n-1}}\sum_{2^{n-1} \pm} \biggl\{\log(\e_1^\pm + \cdots +\e_j^+ + \cdots + \e_n^\pm) + \log(\e_1^\pm + \cdots +\e_j^- + \cdots + \e_n^\pm) \biggr\}  \\
& = \ {1\over 2^{n-1}}\sum_{2^{n} \pm} \biggl\{\log(\e_1^\pm + \cdots +\e_j^\pm + \cdots + \e_n^\pm) \biggr\}   \ =\ 
{1\over 2^{n-1}} \log F(A)          \\
& \geq \ \log(\e_j^-) + \log(\e_j^-).
\end{aligned}
$$ 
Summing over $j=1, ... , n$ gives
$$
{n\over 2^{n-1}} F(A) \ \geq\ \sum_{j=1}^n \bigl \{ \log(\e_j^-) + \log(\e_j^-) \bigr \} \ = \ \log(\det(A))
$$
\qed

\def\e{\epsilon}

\vskip .3in

\noindent
{\headfont  9. Being a G\aa rding-Dirichlet Operator is not Enough for Determinant  Majorization.} 
 
It is not true that the determinant majorization inequality holds for all G-D polynomials.
The Central Ray Hypothesis is crucial.
For a counterexample we take $F:\Sym(\bbr^2) \to \bbr$ to be the cubic polynomial
$$
F(A) \ \equiv\ a_{11}^2 a_{22}.
$$
The polynomial $F$  is $I$-hyperbolic with 
G\aa rading eigenvalues $\l_1^F(A)=\l_2^F(A)=a_{11}$ and $\l_3^F(A)=a_{22}$
since 
$$
F(tI+A) \ =\ (t+a_{11})^2(t+a_{22}).
$$
 The G\aa rding cone $\G$ equals $\{A : a_{11}>0 \ {\rm and}\ a_{22}>0\}$, and
since its closure $\overline \G$ contains $\cp$, this proves that  $F$ is a G-D operator.  
  
  Suppose now that $A>0$ is diagonal,  so that, $a_{12}=0$.   Then
 $$
 {F(A)^{{1\over 3}}\over \det(A)^{{1\over 2}}} \ =\  \left( {a_{11}\over a_{22}}\right)^{{1\over 6}} 
 \ \not \geq \ \g \qquad \text{for any $\g>0$}.
 $$
Of course the Central Ray Hypothesis is not satisfied by the diagonal operator $F(a_{11}, 0, a_{22})$,
so the Basic Lemma 2.1 does not apply.  

We note that since $F$  is a G-D operator, the  inhomogeneous Dirichlet problem
$$
F(D^2 u) \ =\ \left( {\partial^2 u\over \partial x^2}\right)^2  \left( {\partial^2 u\over \partial y^2}\right)  = f(x,y),
\qquad {\partial^2 u\over \partial x^2}, \ {\partial^2 u\over \partial y^2} \geq 0, \ u\bigr|_{\partial\O} = \vf,
$$
with $f\in C(\overline\O), f\geq 0$, and $\phi \in C(\partial\O)$, 
can be uniquely solved on domains $\O\ss\ss\bbr^2$ with smooth strictly convex boundaries
(see [HL$_6$]).

This simple example extends to a family where the determinant majorization fails.

Consider the diagonal G-D operator  $$F(A) \equiv  a_{11}^{N-1} \cdot {1\over n}(a_{22}+ \cdots + a_{n+1, n+1}),$$ 
with $N\geq 3, n\geq 1$, on $\bbr\times \rn$ with coordinates $(t,x)$.  
The G\aa rding cone $\G$ equals $\{a_{11}>0$ and $a_{22}+\cdots + a_{n+1, n+1} >0\}$.
Restricting to diagonal $A>0$ and fixing
$a_{22}, ... , a_{n+1, n+1}$, one has
$$
{F(A)^{1\over N} \over (\det\, A)^{1\over  n+1}} \ =\ a_{11}^{1 - {1\over N} - {1\over n+1}} C \qquad {\rm with}\ \ C>0.
$$
Since $1-{1\over N}-{1\over n+1} \geq 1-{1\over 3} - {1\over 2} = {1\over 6} > 0$, there is no lower bound for the left-hand
term on $A>0$, i.e., determinant majorization fails.

Although, again by [HL$_6$], the continuous inhomogeneous Dirichlet Problem can be solved as above.

When $n\geq 2$, there exists a solution to $F(D^2 u) = k >0$ which is ${2\over N}$-H\"older continuous but no better.
Note that the domain of the operator $F$ is the space of u.s.c.\ functions $u(t,x)$ which are convex in $t$ and $\D$
subharmonic in $x$.

\noindent
{\bf Lemma 9.1. ($n\geq 2$).} {\sl
The singular Pogorelov function $u(t,x) \equiv g(t) |x|^{2\over N}$, where $g$ satisfies $g''(t)^{N-1} g(t) =1$, $g(0)=1, g'(0)=0$, 
is a viscosity solution to}
$$
F(D^2u) \ \equiv\ u_{tt}^{N-1} {1\over n} \D_x u \ =\ k\ \equiv \ {2\over nN} \left (n-2+{2\over N}\right)  \ >0.
\eqno{(9.1)}
$$

\noindent
{\bf Outline of Proof.}  First, if $x\neq 0$, then using ${\partial \over \partial x_j} |x| = {x_j \over |x|}$
one directly computes that ${1\over n} \D_x |x|^{2\over N} = {2\over nN} (n-2+2/N)|x|^{{2\over N}-2}$. 
Hence, $$F(D^2 u) = g''(t)^{N-1} g(t) {1\over n} \D_x  |x|^{2\over N} \ =\ k.$$

Second, use the smooth approximations 
$$
u_\epsilon(t,x) \equiv g(t)(|x|^2+\epsilon)^{1\over N}
\eqno{(9.2)}
$$
which decrease to $u$ as $\e\searrow 0$.

A similar direct calculation yields
$$
F(D^2 u_\e) \ =\ \left({2\over nN}\right) {C|x|^2 + \e n  \over |x|^2 + \e}, \ \ {\rm with}\ \ C \ \equiv \ n-2+{2\over N}.
\eqno{(9.3)}
$$
Note that $k={2\over nN}C$ and that 
$$
{C|x|^2 + \e n  \over |x|^2 + \e} - C \ =\  {\e(n-C)  \over |x|^2+\e},\ \ {\rm with}\ \ n-C \ \equiv \ 2-{2\over N} \geq 1.
$$
Hence, each $u_\e$ is a subsolution of $F(D^2 u)=k$, proving that the decreasing limit $u$ is also a subsolution.

Finally, with $\e>0$ sufficiently small, one can show that there exist $\eta_\e>0$ with $\eta_\e \to 0$ and 
$F(D^2(u_\e - \eta_\e \half |x|^2)) \leq k$ if $u_\e  - \eta_\e \half |x|^2$  is admissible at $x$, i.e., $\D_x u_\e-\eta_\e\geq 0$.
That is, $u_\e - \eta_\e \half |x|^2$ is a supersolution of the equation $F(D^2 u)=k$.  As $\e\to 0$, this converges uniformly to $u$.
So $u$ is also supersolution and therefore a viscosity solution.  \qed

\vskip .3in

\noindent
{\headfont  10. Other Applications of the Basic Lemma.} 

There are other families of polynomial operators on $\Symn$, besides those discussed above, to which the 
Basic Lemma 2.1 does  apply.

Throughout this section we assume that $p(x)$ is a polynomial satisfying the Hypotheses of  Basic Lemma 2.1
and that $E\ss \Symn$ is a closed $\bbr_{\geq0}^n$ monotone set containing $\bbr_{\geq0}^n$.
We also assume the pair $p,E$ satisfies
$$
p(x+y) \geq p(x) \quad \forall\, x\in E\ {\rm  and}\ y\in \bbr_{\geq0}^n.
\eqno{(10.1)}
$$
Since $p$ has coefficients $\geq 0$ (Assumption (1)), $p$ always satisfies (10.1) if we set $E\equiv \bbr_{\geq0}^n$.



\medskip
\noindent
{\bf FAMILY 1. Diagonal Operators.}  Let $F$ be a homogeneous polynomial of degree $N$ on $\Symn$ such that 
$$
F(A) \ =\ p(a_{11}, a_{22}, ... , a_{nn})  \qquad {\rm if}\ \ (a_{11}, ... , a_{nn}) \in E.
$$
where $p$ satisfies the conditions of Basic Lemma 2.1. Then $F$ satisfies 
$$
F(A)^{1\over N} \ \geq \  \det(A)^{1\over n} \qquad \text{for all $A>0$.}
\eqno{(10.2)}
$$

Note that $F$ with domain restricted to $\{A\in\Symn : (a_{11}, ... , a_{nn}) \in E\}$ is a elliptic operator by (10.1).

\pf
By the Basic Lemma 2.1,  $F(A)^{1\over N} \ \geq  (a_{11} a_{22} \cdots a_{nn})^{1\over n}$ for $A>0$, and for 
$A>0$ we have  that $a_{11} a_{22} \cdots a_{nn} \geq \det(A)$.
This last (well known) inequality is proved as follows.  Take the Cholesky decomposition $A=L L^t$ where $L$ is lower triangular.  
Then $\det\, A = (\det\, L)(\det\, L^t) = (\det\, L)^2 = \a_1^2 \cdots \a_n^2$ where $\a_1, ... , \a_n$ are the diagonal entries of $L$.
In terms of row vectors we write $L$ as 
$$
L\ =\ \left(  
\begin{matrix}
w_1 \\
\vdots\\
w_n
\end{matrix}
\right)
     \quad {\rm where}\qquad
 \begin{aligned}
 w_1 \ &= \ (\a_1, 0, 0, ... , 0) \\
  w_2 \ &= \ (*, \a_2, 0, ... , 0) \\
  \vdots  \\
  w_n \ &= \ (*, *,  ... , *, \a_n). \\
 \end{aligned}
$$
Hence, 
$$
A\ =\ \left(  
\begin{matrix}
|w_1|^2 & \  & \ & \ \\
\ & |w_2|^2 & \ &\ * &  \ \\
\  \\
\ & *& \ & \ \\
\ & \ & \ & \ &\  |w_n|^2
\end{matrix}
\right)
$$
and the diagonal entries are $a_{11} \cdots a_{nn} = |w_1|^2 |w_2|^2 \cdots |w_n|^2 \geq \a_1^2 \cdots \a_n^2$.\qed

\medskip

\noindent
{\bf FAMILY 2. Ordered Eigenvalue Operators.}  Let $p(\l_1, ... , \l_n)$ satisfy the assumptions of Basic Lemma 2.1,  and for any $A\in \Symn$ set
$$
F(A) \ \equiv p(\l_1(A),  ... , \l_n(A))  \qquad \text{where $\l_1(A) \leq \l_2(A) \leq \cdots \leq \l_n(A)$.}
$$
Then $F$ satisfies (10.2).  (Note that these operators  $F$ are continuous but generally not polynomials.)

There are many non-symmetric polynomials satisfying these assumptions.

\def\l{x}

\noindent
{\bf Lemma 10.1.}   {\sl
Supoose $q(y_1, ... , y_n)$ and $r(z_1, ... , z_m)$ both satisfy the hypotheses of Basic Lemma 2.1, 
with the same value of $k$ in Assumption (3),  then so 
also does the polynomial}
$$
p(\l_1, ... , \l_{n+m}) \ \equiv\ q(\l_1, ... , \l_n) \, r(\l_{n+1}, ... , \l_{n+m}).
$$

\pf
Note that $p(1, 1 , ... , 1) = q(1, ... , 1) \, r(1, ... , 1) \ =\ 1$.  Also we have
$$
{\partial p \over \partial x_j}(1, ... , 1) \ =\ 
\begin{cases}
{\partial q \over \partial x_j}(1, ... , 1)\, r(1, ... , 1) \ =\ k \qquad{\rm if}\ \ 1\leq j\leq n  \\
q(1, ... , 1)\, {\partial r \over \partial x_j}(1, ... , 1)  \ =\ k\qquad{\rm if}\ \ n+1\leq j\leq n  +m.
\end{cases}
$$
It is clear that the coefficients of $p$ are all $\geq 0$.\qed

\noindent
{\bf Lemma 10.2.}   {\sl
Supoose $q(y_1, ... , y_n)$ and $r(z_1, ... , z_m)$ both satisfy the hypotheses of Basic Lemma 2.1 
and have the same degree.  Let $k$ and $k'$ be the constants in Assumption (3) for $q$ and $r$ respectively.
The the polynomial
$$
p(\l_1, ... , \l_{n+m}) \ \equiv {k'\over k+k'} q(\l_1, ... , \l_n)  +  {k\over k+k'} r(\l_{n+1}, ... , \l_{n+m}).
$$
also satisfies  the hypotheses of Basic Lemma 2.1.}

\pf
We clearly have that $p$ is homogeneous, and $p(1, ... , 1)=1$. If $j\leq n$, 
$$
{\partial p \over \partial x_j}(1, ... , 1) \ =\ {k'\over k+k'} {\partial q \over \partial x_j}(1, ... , 1)\ =\  {k'k\over k+k'},
$$
and if $j\geq n+1$, then
$$
{\partial p \over \partial x_j}(1, ... , 1) \ =\ {k\over k+k'} {\partial r \over \partial x_j}(1, ... , 1)\ =\  {k'k\over k+k'}.
$$
Again it is clear that the coefficients of $p$ are all $\geq 0$.\qed

\def\l{\lambda}

\noindent
{\bf Example 10.3.}  The polynomials  $P_j(x) = {1 \over j}\{x_1^j + \cdots + x_j^j\}$  for $j\geq 1$
each  satisfy the hypotheses of Basic Lemma 2.1, 
with the same value $k=1$ in Assumption (3). Hence Lemma 10.1 can be applied to any two
of these, and in fact to any finite number taken from this set. Consider the simplest case.

For $A\in\Sym(\bbr^3)$ with eigenvalues $\l_1(A)\leq \l_2(A)\leq \l_3(A)$, and $p(\l) \equiv P_1(\l_1) P_2(\l_2, \l_3)$
$$
F(A) \ =\ P_1(\l_1(A)) P_2(\l_2(A), \l_3(A)) \ =\ \l_1(A)\half \{\l_2(A)^2 + \l_3(A)^2\}.
$$
The connected component of the set $\{A : F(A)>0\}$ containing $I$ is $\Int \cp\equiv \{A: A>0\}$.  On this set $F>0$, and $F=0$  
on $\partial \cp$.  The eigenvalues $\l_2(A), \l_3(A)$ are positive and $\l_2(A)^2 + \l_3(A)^2\geq 2 \l_1(A)^2$ on this set.
 Here the determinant majorization is clear.

\vskip .5in

\def\e{\epsilon}
\noindent
{\headfont  Appendix A. The Interior of the Polar Cone.}    

Here we describe some useful criteria for determining when a vector is in the interior of the polar of
a convex cone.  We start with a non-empty open convex cone  $\G\ss V$ in a finite dimensional 
real vector space $V$.  The polar
$$
\G^0\ \equiv \  \{y\in V^* : \bra yx \geq 0 \ \forall \, x\in \G\}
$$
is a non-empty closed convex cone.  As such it has interior when considered as a cone in its vector space
span $S\ss V^*$.  We denote this relative interior by $\G^*$ and shall refer to it as the {\bf open polar}, keeping in mind
that it is an open convex cone  in its span $S$, which may be have  lower dimension than that of $V$.
The orthogonal complement $E\equiv S^\perp$ is called the {\bf edge}  of  $\G$.  
It is characterized as the linear subspace of $\G$ which contains all the lines (through the origin) in $\G$.
Note that if $\overline \G$ is self-polar (i.e., $\G^0 = \overline \G^0 = \overline \G)$, then $\G^*=\G$.  In particular, $(\Int \cp)^* =\Int \cp$.

The simple 2-dimensional example $\G\equiv \{x\in \bbr^2 : x_1>0, x_2 >0\}$, where $E=\{0\}$
and $\G^*=\G$ (self-polar), is a counterexample to the statement:
$$
{\rm For} \ \ y\in S:  \qquad \ \bra yx >0\ \ \forall\, x\in \G \quad \implies \quad y\in \G^*,
$$
(take $y=(1,0)$.)
However, various strengthenings of this display provide useful criteria for $y$ to belong to the open polar $\G^*$.

\noindent
{\bf Lemma A.1. (The Open Polar Criteria).}    {\sl
Suppose $y\in S$.  The following are equivalent.

(1) \ \ $y \in \G^*$.

(2) \ \ $\exists\,  \e>0 \ \ \text{such that } \ \ \bra yx \geq \e|x|\ \ \forall\, x \in \overline\G \cap S$.

(3) \ \ $\bra yx >0 \quad \forall \, x\in \overline \G - E$.

(4) \ \ $\bra yx >0\quad \forall\, x\in (\partial \G-\{0\})\cap S$.

(5) \ \  $\bra yx >0\quad \forall\, x\in \overline \G \cap S, \ x\neq 0$.
}

 \noindent
 {\bf Proof.} (1) $\iff$ (2):  Note that:  
 $$
 \begin{aligned}
 y\in \G^* \ &\iff \ \text{there exist an $\e$-ball $B_\e(y) \ss \G^0$ about $y$} \\
 &\iff \exists \e>0 \ {\rm such\  that}\ \bra{y+\e z}x \geq 0 \ \forall\,|z|\leq 1, x\in \G \\
&\iff  \text{(a): $\exists \, \e>0 \  {\rm such\  that}\ \bra yx \geq \e \bra zx \ \ \forall\, |z| \leq 1, x\in\G$.}
  \end{aligned}
$$ 
Taking $z=x/|x|$ in (a) yields (2),  while (2) yields (a) since $\bra zx \leq |x|$ if $ |z|\leq 1$.

To see  (2) $\Rightarrow$ (3), suppose $x\in \overline\G - E$ and decompose $x$ into $x=a+b$ where
$a\in E, b\in S$.  Then $b\in \overline \G\cap S$ and $b\neq 0$, so $\bra yx = \bra bx >0$, which proves (3).

That (3) $\Rightarrow$ (4) is obvious.

Since $\overline G$ is convex, it is easy to see that (4) $\Rightarrow$ (5).

To see  (5) $\Rightarrow$ (2): We can  assume that $E=\{0\}$ and $S=V$.  Obviously, (5) $\Rightarrow$ $y\in \G^0$.  Hence,
$$
\e \ \equiv \ \inf_{x\in \overline \G, |x|=1} \bra yx \ \geq\ 0.
$$ 
Since $\overline \G \cap \{|x|=1\}$ is compact, there exists $x_0 \in \overline \G \cap \{|x|=1\}$ with $\e= \bra y{x_0}$.
Now (5) $\Rightarrow$ $\e>0$.   Hence, $\bra y {{x\over |x|}} \geq \e \ \ \forall\, x\in \overline\G-\{0\}$, or 
$\bra yx \geq \e|x|\ \ \forall\, x\in \overline\G$.  \qed

\vskip .3in
\noindent
{\headfont  Appendix B. The G\aa rding Gradient Map.} Let $F$ be a G-D polynomial of degree $N$ on $\Symn$, 
and let $\G$ be its G\aa rding cone.  We fix $A\in \G$ and write (by Elementary Property (3) in [HL$_3$])
$$
F(A+tB) \ =\ F(A)\prod_{j=1}^N(1+t\l^{F,A}_j(B))
$$
where $\l_j^{F,A}(B)$ are the $F$-eigenvalues of $B\in \Symn$ with respect to $A\in \G$. Hence,
$$
{1\over F(A)} {d \over d t}  \log \, F(A+tB) \biggr|_{t=0} 
\ =\ 
\sum_{j=0}^N \l_j^{F,A}(B) \ =\  {1\over F(A)} \bra{(\nabla F)_A}{B}
\eqno{(B.1)}
$$

We shall assume that $F$ is {\bf complete} which means that all the variables in $\rn$ are needed to define the G-D operator $F$.
There are many useful equivalent ways of describing completeness (see section 3 of [HL$_5$] for a full discussion of completeness).
We point out that an invariant G-D operator is complete.

Here we are using the notation $(\nabla F)_A \equiv D_A F$.

\noindent
{\bf Proposition B.1.}  {\sl A complete G-D operator $F$ has uniformly  elliptic linearization at each point of its G\aa rding cone $\G$.
That is, at $A\in\G$
$$
B \mapsto \text{$\bra {(\nabla  F)_A}{B}$ has coefficient matrix $(\nabla F)_A > 0$.}
$$
}
 \pf
 We use (3) in Proposition 3.5 in [HL$_5$]  saying that completeness is equivalent to 
 $$
 \cp\cap E \ =\ \{0\}.
 \eqno{(B.2)}
 $$
Note that
$$
B>0 \quad \iff\quad \bra BP \ >\ 0\ \ \forall P\in \cp-\{0\}.
\eqno{(B.3)}
$$
Now (B.1) for the derivative of $F$ says that for $A\in \G$ and $P \in \Int \cp\ss \G$
$$
\bra {(\nabla F)_A}{P}  \ =\ F(A) \sum_{j=1}^N \l_j^{F,A}(P).
$$
Now $F(A)>0$  and $\l_j^{F,A}(P) \geq 0$ for all $P\geq 0$ since $\cp\ss\overline \G$.
Hence, $\bra {(\nabla F)_A}{P} \geq 0$, and $=0 \ \iff \ \l_1^{F,A}(P) = \cdots = \l_N^{F,A}(P) = 0$, but
$\l_1^{F,A}(P) = \cdots = \l_N^{F,A}(P) = 0 \ \iff\ P \in E$.  Since $P\geq 0$ and $P\in E$ the completeness hypothesis (B.2) implies $P=0$.  Applying (B.3) we have $(\nabla F)_A > 0$.
\qed

\vskip .3in

\def\gg{{\mathfrak g}}
\noindent
{\headfont  Appendix C. The Exhaustion Lemma.} 

Again we take the general point of view.
Let $\G\ss V$ be the G\aa rding cone for  a G\aa rding polynomial $\gg$ on a finite dimensional vector space $V$.
It is useful to construct convex {\sl exhaustion functions} $\psi(x)$ for $\G$,  that is, functions $\psi\in C^\infty(\G)$
which are convex and all the prelevel sets $K_c \equiv \{x\in\G : \psi(x)\leq c\}$ are compact.
Now $\psi(x) \equiv - \log\, \gg(x) \in C^\infty(\G)$ is convex and extends to a continuous function on $\overline \G$
which is equal to $+\infty$ on $\partial \G$.  However, the associated sets $K_c = \{x\in \G : \gg(x) \geq e^{-c}\}$ are not compact.
For example, when the dimension $n=1$, if $\gg(x) =x, \G=(0, \infty)\ss V\equiv \bbr$,
these prelevel sets are not compact.  Modifying $-\log\, x$ by adding the linear function $\g x$ with $\g>0$ a constant, 
yields an exhaustion function $\psi$ for $\G$.

There is a caveat here.   This scenario given above, of adding a linear function to $\log\, \gg$ does not work if $\G$ has an edge $E$.
Recall from Appendix A that $E$ is the largest linear subspace contained in $\G$.  For an extreme example, take the Laplacian
$\gg(A)= \tr(A)$ on $\Symn$ where $\G$ is the half-space $\{ \tr(A) \geq 0\}$ and $E=\{\tr(A)=0\}$.   If we divide by the edge, we fall onto the example above and the argument does work there.  In general, the  proposed compactness argument 
 holds only after dividing by the edge, or equivalently, restricting to the span $S=E^\perp$. 
 We elect to restrict  $\G$ to the span $S$ since this conforms with the common choice for the complex Monge-Amp\`ere
 operator, where $S = \Sym_\bbc(\bbc^n)$ is the Hermitian symmetric matrices.

\noindent
{\bf Definition C.1.} The convex cone $\G\ss V$    is called {\bf regular} if the edge of $\overline\G$ is $\{0\}$.

\noindent
{\bf The Exhaustion Theorem C.2.} {\sl
Suppose $\gg$ is a G\aa rding polynomial on $V$ of degree $N$ with G\aa rding cone $\G$.
Fix $y \in \G^*$, the (relatively) open polar.  First, assume $\G$ is regular.
Then the function $\psi(x) \equiv \bra yx - \log\, \gg(x)$ is a strictly convex exhaustion function for $\G$,
so the prelevel sets
$$
K_c \ \equiv\ \{x\in \G : \bra yx -\log \, \gg(x) \leq c\},\qquad\text{for $c\in \bbr$}
$$
are compact.

If $\G$ is not regular, then restricting $\psi$ to the span $S\equiv E^\perp$ gives a strictly convex exhaustion function
for the regular cone $\G_S \equiv \G\cap S$.  Hence, the prelevel sets
$$
K_c \ \equiv\ \{x\in \G\cap S: \bra yx -\log \, \gg(x) \leq c\},\qquad\text{for $c\in \bbr$}
$$
are compact subsets of $\G_S = \G\cap S$.
}

\pf
It suffices to prove the theorem when $\G$ is regular, since otherwise $\gg\bigr|_{S}$ has regular G\aa rding cone $\G_S = 
\G\cap S$ with the same open polar $\G^*_S = \G^*$.

Recall that $-\log\, \gg(x)$ has second derivative at a point $x\in \G$ given by
$$
\left\{  D_x^2(-\log \, \gg)   \right\} (\xi, \xi) \ =\ \sum_{j=1}^N\left ( \l_j^{\gg, x}(\x)   \right)^2
$$
for all $\x \in V$, where the G\aa rding eigenvalues are taken with respect to the direction $x\in \G$.
By G\aa rding [G\aa r] (or see  [HL$_4$]), the nullity set
$$
\{ \x \in V : \l_1^{\gg, x}(\x) = \cdots + \l_N^{\gg, x}(\x) = 0\}
$$
 equals the edge $E$.
Hence,  the function $-\log\, \gg(x)$ is strictly convex on $S=E^\perp$, and so $\psi(x)  \equiv \bra y x -\log\, \gg(x)$ 
has the same property  since 
$\bra y x$ is affine.  

Notice that since $\gg \equiv 0$ on $\partial \G$  and $\bra yx$ is finite,
the function $\psi \equiv +\infty$ on $\partial \G$.
We conclude that $K_c$ is a closed  subset of $\G$.

It remains to  show that $K_c$ is bounded.  For this we use the full hypothesis that $y$ belongs to the open polar  $\G^*$
of $\G$, which equals the interior of $\overline\G^0$.
By Lemma A.1 (2)  this is equivalent to the statement

\noindent (2) \ \ $\exists\ \e>0$  such that $\bra yx  \geq \e|x|\ \forall\,x\in \overline \G$.

This implies that:
$$
\begin{aligned}
K_c \ &\ss\ \{x\in  \overline \G : e^{\e|x|} \leq e^c \gg(x)\} \ss \left \{x\in  \overline \G : {(\e|x|)^{N+1} \over (N+1)!} \leq e^c\gg(x) \right \}  \\
&= \left \{x\in  \overline  \G : |x| \leq {(N+1)! e^c \over \e^{N+1}} \gg\left({x\over |x|}  \right)\right\} \ \ss\ B_R(0)
\end{aligned}
$$
with
$$\qquad
R \equiv {(N+1)! e^c\over \e^{N+1}} \sup_{\x\in \overline\G, |\x|=1} \gg(\x)\ <\ \infty.    \hskip.5in \mathqed
$$
 
 \noindent
 {\bf Corollary C.3.} {\sl
 For constants $c_1, c_2\in \bbr$ and $y\in\G^*$, the set $$K_{c_1, c_2}  \equiv \{x\in \G : \bra y x \leq c_1 \ {\rm and}\ \log\, \gg(x) \geq c_2\}$$
 is compact.  
 }

\pf One has $K_{c_1, c_2} \ss K_{c_1 - c_2}$.\qed

\vskip .3in

\noindent
{\headfont  Appendix D. The G\"uler  Derivative Estimates for G-D Operators.} 

\smallskip

Suppose $F$ is a G\aa rding-Dirichlet operator on $\Symn$ of degree $N$ with   G\aa rding cone $\G$.
(We assume $F\not \equiv 0$, and recall that  by definition $\G$ is open.) For each $A\in \G$ there is a simple formula for the $k^{\rm th}$
derivative of  $L\equiv \log\, F$ at $A$ in terms of the G\aa rding eigenvalues $\l_j^{F,A}$ of $F$ w.r.t.\  the direction $A$.

In this section we shall abbreviate $\l_j^{F,A}$ to $\l_j^{A}$ or just $\l_j$.

\noindent
{\bf Theorem D.1. (Derivatives).} {\sl
At a point $A\in\G$ and $\forall \,B\in \Symn$:

 (1)
$$
\begin{aligned} 
\bra {D_AF} B \ &=\ F(A) \sum_{j=1}^N \l_j^A(B) \qquad\qquad{\rm and} \\
\bra {D_A F} B\ \geq\ 0 \ \ \forall\, B\in\overline \G, & \ \qquad \text{with equality $\iff B\in E$.}
\end{aligned}
$$
\indent (2)
 $$
 \begin{aligned} 
D_A^2(\log\, F)(B, B) \ =\  -  \sum_{j=1}^N&( \l_j^A(B))^2 \qquad\qquad{\rm and} \qquad\  \\
\text{the quadratic form } \quad D_A^2(\log\, F) \ \leq\ 0\  & \ \qquad \text{with null space  $E$.}
\end{aligned}
 $$
\indent(3)
 $$
 \begin{aligned} 
D_A^3(\log\, F)(B, B,B) \ =\  2  \sum_{j=1}^N&( \l_j^A(B))^3  \quad
\end{aligned}
 $$
\indent(4) 
 $$
 {{
D_A^{(k)}(\log\, F)(B, ... ,B) \ =\  (-1)^{k-1} (k-1)! \sum_{j=1}^N( \l_j^A(B))^k \quad \forall\, k\geq 1
}}
 $$
 }
 
In particular,
\begin{itemize}
\item [(1)$'$] 
The directional derivative of $F$ at $A$ is strictly increasing in all directions $B$ in the closed G\aa rding cone
$\overline\G$, other than the edge directions where it is zero.
\item[(2)$'$]
The restriction of $L = \log\, F$ to the span $S$ is a strictly concave operator on $\G_S \equiv S\cap \G$.
\end{itemize}

\pf
For $A \in \G$ and $B\in \Symn$, let $\vf(t) \equiv \log\, F(A+tB)$.  Then
$$
F(A+tB) \ =\  F(A) \prod_{j=1}^N (1+t \l^A_j(B)).
$$
(See for example  [G\aa r] or  the  Elementary Property (3) in  [HL$_3$].)  

Since $\vf(t) = \log\, F(A) + \sum_{j=1}^N \log\, (1+t\l_j^A)$, we have (with $\l_j\equiv \l_j^A$) that
$$
\vf'(t) \ =\ \sum_{j=1}^N {\l_j \over 1+t\l_j}, \ \ \vf''(t) \ =\ - \sum_{j=1}^N {\l_j^2 \over (1+t\l_j)^2},
$$
and for $k\geq 2$, 
$$
\vf^{(k)}(t) \ =\ (-1)^{k-1} (k-1)! \sum_{j=1}^N{\l_j^k\over (1+t\l_j)^k}.
$$
Since $\vf^{(k)}(0) = \{D^{(k)}_A(\log\, F)\}(B, ... , B)$ this proves all the equalities in (1), (2), (3) and (4).
The remainder of the theorem follows from the fact that for all $B\in\overline \G$, the eigenvalues 
 $\l_1^A(B), ... , \l^A_N(B)$ are all $\geq 0$ 
with equality iff $B\in E$,  that is , the nullity $\{B : \l_1^A(B) = \cdots \l_N^A(B) = 0\}$ equals the edge $E$
(see [G\aa r] or [HL$_3$]).\qed

Now with $\l_j \equiv\l_j^A(B) \equiv \l_j^{F,A}(B)$ and $\l \equiv (\l_1, ... , \l_N)\in \bbr^N$, the little $\ell^k$-norm of $\l$ is 
$$
\|\l\|_k  \ \equiv \ \left( \sum_{j=1}^N  \| \l_j\|^k  \right)^{1\over k} \ \ \  {\rm for} \ k\geq 1 \qquad{\rm and} \qquad \| \l \|_\infty = \sup_{1\leq j\leq N} |\l_j|.
$$
Recall that 
$$
\|\l \|_k \ \leq \| \l \|_\ell \quad {\rm for}\ \ \ell=1, 2, ... , k-1\ \ \text{with equality iff $\l$ is an axis vector}.
\eqno{(D.1)}
$$
To see this  we can assume that $\l\neq 0$.  By homogeneity one can assume $\|\l\|_\ell=1$, and therefore $0\leq |\l_j| \leq 1$
 for all $j$.  Then $\ell< k$ implies $\sum |\l_j|^k \leq \sum|\l_j|^\ell = 1$. Hence, $\|\l\|_k \leq 1 = \|\l\|_\ell$ with equality
if and only if $\l$ is a unit axis vector.
 
 This translates into an upper bound estimate for the $k^{\rm th}$ derivative by lower order derivatives.
 
\noindent
{\bf Theorem D.2. (G\"uler [G\"ul]).}  {\sl
For $A \in \G$, $B \in \Symn$ and $k$ fixed,
if $\ell = 2, 4, ...$ is even and $\ell \leq k$, then
$$
 \left|   {1\over (k-1)!} \left \{ D^{(k)}_A  \,\log\, F \right\} (B, ... , B)   \right|^{1\over k}
\ \leq \ 
\left[  {1\over (\ell-1)!} \left  \{ D^{(\ell)}_A \, \log\, F \right \} (B, ... , B)   \right]^{1\over \ell}
$$
}

\pf
Note that for $k=1, 3, ...$ odd,
$$
 \left|   {1\over (k-1)!} \left  \{ D^{(k)}_A \, \log\, F \right \} (B, ... , B)   \right| = \left|   \sum_{j=1}^N \l^k_j  \right| 
  \leq \sum_{j=1}^N | \l_j|^k \ =\   \|\l\|^k_k,
$$
and for $k=2, 4, ...$ even 
$$
 {1\over (k-1)!} \left  \{ D^{(k)}_A \log\, F  \right \} (B, ... , B) \ =\ \|\l\|_k^k.
$$
Hence, in both cases,  $\|\l\|_k \leq \|\l\|_\ell$ for $\ell=1, ... , k-1$ yields the estimates for the $k^{\rm th}$ derivative of $\log\,F$.
\qed

For equality to hold in Theorem D.2, $B$ must either be in the edge $E$ or have exactly one non-zero $A$-eigenvalue.

\medskip

\centerline{\bf Further Discussion --  Second Derivative Formulas}  

The above first and second derivative formulas for $\log\, F$ at $A\in \G$, namely
$$
(1) \ \ \bra{D_A \log\, F}{B}\ =\ \sum_{j=1}^N \l_j
\and
(2) \ \ (D_A^2 \log\,F)(B,B) \ =\ - \sum_{j=1}^N \l_j^2
$$
for all $B\in \Symn$ can be used to conclude interesting formulas for the 
second derivatives of $F$ and $F^{1\over N}$.
Here we abbeiviate $\l_j^{F,A}(B)$ to just $\l_j$ or $\l_j(B)$, and we recall the discriminant polynomial
$${\rm Discr}(\l_1, ... ,  \l_n) = \sum_{i<j}(\l_i-\l_j)^2.$$
\noindent
{\bf Proposition D.3.}

(D.1) \ \ \ $(D_A^2 \log\, F)(B,B) \ =\ -|\l(B)|^2$

(D.2) \ \ \ $(D_A^2 F)(B,B) \ =\ 2F(A) \s_2(\l(B))$

(D.3) \ \ \ $(D_A F^{1\over N})(B,B) \ =\ -{1\over N^2} F(A)^{1\over N} {\rm Discr}(\l(B))$

{\sl
Moreover, the quadratic forms $D^2_A \log\, F$ and $D^2_A F$ have the same null space, namely
$N \equiv \{B : \l_1(B)= \cdots =\l_N(B)=0\}$.  The nullity $N$ is the same as the edge $E\equiv \overline{\G}\cap (-\overline {\G})$
and also the linearity of $F$ which is, by definition, the largest linear subspace on which $F=0$.
Modulo this subspace, $D_A^2 \log\, F <0$, while $D_A^2 F$  has Lorentzian signature with future light cone $\G$.
Fiinally, the quadratic form  Discr$(\l(B))\geq 0$ with null space  $\{B \in\Symn : \l_1(B)=\cdots = \l_N(B)\}  = E+\bbr\cdot A$
of dimension = dim$(E)+1$.}

Note that (D.3) proves the G\aa rding Lemma that $F^{1\over N}$ is concave.

\noindent
{\bf Proof of (D.2):}  By (1) and (2) above,
$$
\begin{aligned}
-\sum_{j=1}^N \l_j^2 \ &=\ D_A^2 \log\,F \ =\ D_A\left( {D_AF\over F(A)}\right) \ = \ {D_A^2(F)\over F(A)} \\
&= \  - {D_AF \circ D_AF \over F(A)^2}\ =\ {D_A^2F\over F(A)} - \left(\sum_{j=1}^N \l_j\right)^2.
\end{aligned}
$$
Note that $(\sum_j\l_j)^2 - \sum_j\l_j^2 = 2\s_2(\l)$. \qed

\noindent
{\bf Sketch of proof of (D.3):}  
Use the standard formula (obtained by expanding out the right hand term)
$$
N |\l(B)|^2  - \s_1(\l(B))^2 \ =\ {\rm Descr}(\l(B)) \ \equiv \  \sum_{i<j} (\l_i(B) - \l_j(B))^2
$$
to compute that 
$$
(D_A^2 F^{1\over N}) (B, B) \ =\ - {1\over N^2}  F^{1\over N} (A)   \left\{     {\rm Descr}(\l^A(B))     \right\}.
$$

\vskip .3in

\noindent
{\headfont  Appendix E.  The Central Ray Hypothesis.} 

\smallskip

\def\cR{{\mathcal R}}

Suppose $F$ is a G-D operator on $\Symn$ of degree $N$ with G\aa rding cone $\G$.
(We assume $F\not \equiv 0$ and by definition $\G$ is open and contains $\{A: A>0\}$.)

The next lemma shows that there exists a unique ray $\cR$ contained in the G\aa rding cone $\G$,
with the property that $\G$ is, in a certain sense, symmetrical about $\cR$.

\noindent
{\bf Lemma E.1.}  
$$
F(A)^{1\over N} \ \leq \left [     \sup_{B\in \overline\G, \|B\|=1}  F(B)^{1\over N}    \right]
\|A\|\qquad \text{for all $A\in \overline \G-E$}.
\eqno{(E.1)}
$$

{\sl 
This sup is attained at a unique maximum point $B_0$ which belongs to $\G_S \equiv \G\cap S$
where $S$ is the span.  Equivalently, equality occurs at a point $A\in\overline\G \iff A\in \cR \equiv \bbr^+\cdot  B_0$, the ray 
through $B_0$.
}

\noindent
{\bf Definition E.2.}  This unique ray $\cR$ for $(F,\G)$ is called the {\bf central ray}.

\pf Suppose $B_0$ is a maximum point.  First, $B_0\in\G$ since $F$ vanishes on $\partial \G$.
Secondly, $B_0\in S$.  If not, then writing  $B_0 = B_E+B_S$ with respect to the orthogonal decomposition
$\Symn = E\oplus S$, we have $B_E \neq 0$.  Therefore, $\|B_S\|< \|B_0\|=1$, and hence
$$
F\left( {B_S \over \|B_S\|}      \right)^{1\over N} 
\ =\ 
 {1 \over \|B_S\|} F(B_S)^{1\over N}
 \ =\   {1 \over \|B_S\|} F(B_0)^{1\over N} 
 \ > \  F(B_0)^{1\over N}.
$$
which contradicts that $B_0$ is a maximum point.  

For the remainder of the proof, we can restrict to the span $S = E^\perp$, and we set 
$\g \equiv    \sup_{B\in \overline\G, \|B\|=1}  F(B)^{1\over N}$.

The tangent affine hyperplane $H\equiv \{A : \bra {B_0} A =1\}$ separates the level sets $\{F(B)^{1\over N} = \g\}$ and $\{\|B\|=1\}$
since $F^{1\over N}$  is strictly concave, and $\|B\|$ is strictly convex on $S$.  Furthermore,  $B_0$ is the only point
of intersection of these hypersurfaces, since $H$ touches each of these hypersurfaces only at the point $B_0$.

If the sup in (E.1) is attained at a point $A$ we may assume $\|A\|=1$ by dividing the inequality by $\|A\|$.
Now we must have $A=B_0$ since all other points on the unit sphere are strictly below the hypersurface $F^{1\over N}=\g$
\qed

\noindent
{\bf Corollary E.3.} {\sl There is a unique point $A_0$, up to positive multiples, where the sphere and the level set of $F^{1\over N}$ 
have a common  normal.}

\pf
Let $A_0$ be a point where there is a common normal.  We can renormalize so that $\|A_0\|=1$, since everything is of degree 1.
  Let $H_0$ be the tangent hyperplane to the sphere at $A_0$.  Then as in the proof of the theorem above $H_0$ strictly separates
the sphere and the hypersurface $F(A)^{1\over N}= F(A_0)^{1\over N}$ away from $A_0$; and the value of 
$F^{1\over N}$ at all points  $\neq A_0$ in the sphere are $<  F(A_0)^{1\over N}$.  Thus if $A_0\neq B_0$ we would have
$  F(A_0)^{1\over N}>  F(B_0)^{1\over N}>   F(A_0)^{1\over N}$. \qed

\noindent
{\bf Corollary E.4.} {\sl The central ray $\cR$ is characterized by the fact that for any non-zero point $B\in\cR$, $D_BF$ and $D_B \|A\|^2 = 2B$ are positively  proportional, i.e., $$D_B F = kB$$ with $B\in \G_S = S\cap \G$ and $k>0$.}

We now examine conditions which imply that the central ray is $\bbr_+ I$.

\noindent
{\bf Propositions E.5.}  {\sl
If $F$ is invariant under O$(n)$, SU$(m)$ with $n=2m$,  Sp$_m$ with $n=4m$, or any other subgroup of O$(n)$ whose only fixed
line in $\Symn$ is the one through the identity $I$, then the central ray $\cR$ is the ray through $I$.
}

\pf  If $F$ is invariant under a subgroup $G\ss {\rm O}(n)$, then $G$ fixes the central ray.\qed

\noindent
{\bf Definition E.6.}  The {\bf Central Ray Hypothesis (CRH)} is that: The identity $I$ generates the central ray $\cR$ of $(\G, F)$.

There are several important ways of formulating the CRH, which we now describe.

\noindent
{\bf Definition E.7.}  The {\bf  G\aa rding  (or $F$-) Laplacian}, denoted $\D^F$, is defined by
$$
\D^F(B) \ \equiv\ \sum_{j=1}^N \l_j^{F,I}(B) \ \equiv \  \tr^F(B).
$$
Note that, after normalizing $F(I)=1$, 
$$
\begin{aligned}
\D^F(B) \ &=\ {d\over dt}\biggr|_{t=0} F(I+tB) \ =\ {1\over (N-1)!} {d^{(N-1)} \over dt^{(N-1)}}\biggr|_{t=0} F(tI+B) \\
&=\ \s_1^F(B) \ =\ \bra {D_I\, \log F}{B}. 
\end{aligned}
$$

\noindent
{\bf Theorem E.8. (Equivalent Formulations of the CRH).}  {\sl    
The following conditions are equivalent.

(1)
$$
\begin{aligned}
F(A)^{1\over N} \ \leq \left[ \sup_{\|B\|^2=n,B\in \G}  F(B)^{1\over N}   \right] & \|A\| \qquad \forall\, A\in \G \\
&\text{with equality iff $A=tI$, $t>0$.}
\end{aligned}
$$

(2)\ \ $D_I\, \log\, F\ =\ {1\over F(I)} D_I F \ =\ kI, \quad {\rm for} \ \ k>0$.

(3) \ \ $\D^F \ =\ k \D^{\rm standard},  \quad {\rm for} \ \ k>0$.

(4)\ \ $F(A)^{1\over N} \ \leq  F(I)^{1\over N} {k\over N} \bra IA\quad  \forall\, A\in \G \qquad \text{with equality iff} \ A\in tI+E$

\noindent
for $t>0$. The $k,n, N$ are related by $kn=N$.
}

 \pf
 We have already shown that (1) $\iff$ (2).  Now (2) $\iff$ (3) since $\D^F(B) = \bra{D_I \log\, F}{B}$ implies
 that $\D^F = k \D^{\rm standard} \iff  D_I \log\, F = kI$. 
  Statement (4) says that $F^{1\over N}(A)$ is bounded above by $c\bra IA$  with equality at $tI+E, t>0$. 
 Thus the derivative of $F^{1\over N}$ at $I$ is a multiple of $I$ which is equivalent to (3).\qed

\vfill\eject

\centerline{\headfont References} 

\noindent
[AO] S. Abja and G. Olive, {\sl Local regularity for concave homogeneous complex
degenerate elliptic equations dominating  the
Monge-Amp\`ere equation},   Ann.\ Mat.\ Pura Appl.\ (4) {\bf201} (2022) no.\ 2, 561-587.

\noindent
[ADO]  S. Abja, S. Dinew and G. Olive, {\sl Uniform estimates for concave
homogeneous complex degenerate elliptic equations comparable to the Monge-Amp\`ere
equation}, Potential Anal. (2022). \newline https://link.springer.com/article/10.1007/s11118-022-10009-w

\noindent
[AV]  S. Alesker and M. Verbitsky, {\sl Plurisubharmonic functions on hypercomplex  manifolds and HKT-geometry},
 ArXiv:0510140v3

\noindent
[BGLS]   H. H.  Bauschke, O.  G\"uler, A. S.  Lewis, H. S.  Sendov, {\sl Hyperbolic polynomials and convex
analysis},  Canad. J. Math. {\bf 53} (2001), no. 3, 470-488.

\noindent
[CC] L.  Caffarelli and X. Cabr\'e, {\sl Fully nonlinear elliptic equations}, American
Math. Soc. Colloquium Publications, vol. 43, American Math.
Soc., Providence, RI, 1995.

\noindent
[CNS] L. Caffarelli, L. Nirenberg, and J. Spruck, {\sl The Dirichlet problem for nonlinear
second-order elliptic equations. III. Functions of the eigenvalues of the Hessian},
Acta Math. {\bf 155} (1985), no. 3-4, 261-301.

\noindent
[D] S. Dinew, {\sl Interior estimates for p-plurisubharmonic functions},
arXiv:2006.12979.

\noindent
[G\aa r] L. G\aa rding, {\sl An inequality for hyperbolic polynomials}, J. Math. Mech. 8 (1959),
957-965.

\noindent
[G\"ul]  O. G\"uler, {\sl Hyperbolic polynomials and interior point methods for convex programming}, 
Mathematics of Operations Research, {\bf 22} No. 2, 1997, 350-377.

\noindent
[Gur]  l. Gurvits, {\sl Van der Waerden/Schrijver-Valiant like conjectures
and stable (aka hyperbolic) homogeneous
polynomials :
one theorem for all},  Electron. J. Combin. {\bf 15} (2008), no. 1, Research Paper 66, 26 pp.

\noindent
[GPT] B. Guo, D.H. Phong, and F. Tong, {\sl On $L^\infty$-estimates for complex Monge-Amp`ere equations},
arXiv:2106.02224.

\noindent
[GP$_1$] B.Guo and D. H.Phong, {\sl On  $L^\infty$-estimates for fully nonlinear partial differential equations on hermitian manifolds}, 
  arXiv:2204.12549v1.

\noindent
[GP$_2$]  \ \------------ , {\sl Uniform entropy and energy bounds for fully non-linear equations}, ArXiv:2207.08983v1, July 2022.

\noindent
 [HL$_1$] F. R. Harvey and H.B. Lawson, {\sl Dirichlet duality and the nonlinear
Dirichlet problem}, Comm. Pure Appl. Math. 62 (2009), no. 3, 396-443.

\noindent
[HL$_2$] \ \------------ , {\sl Dirichlet duality and the nonlinear Dirichlet problem on Riemannian manifolds},
J. Diff. Geom. 88 (2011), no. 3, 395-482.

\noindent
[HL$_3$] \ \------------ , {\sl  Hyperbolic polynomials and the Dirichlet problem},    ArXiv:0912.5220.

\noindent
[HL$_4$] \ \------------ ,   {\sl  G\aa rding's theory of hyperbolic polynomials},
   {\sl Comm. in Pure  App. Math.}  {\bf 66} no. 7 (2013), 1102-1128.
  
\noindent
[HL$_5$] \ \------------ ,  {\sl Lagrangian potential theory and a Lagrangian equation of Monge-Amp\`ere type}, ,pp. 217- 257 in Surveys in Differential Geometry, edited by H.-D. Cao, J. Li, R. Schoen and S.-T. Yau, vol. 22, International Press, Somerville, MA, 2018.

\noindent
[HL$_6$]    \ \------------ ,    {\sl The inhomogeneous Dirichlet problem for natural operators on manifolds}, 
Annales de l'Institut Fourier, {\bf 69} No. 7 (2019),  3017-3064.
ArXiv:1805.11121.

\noindent
[W]  Y. Wang, {\sl On the $C^{2,\a}$-regularity of the complex Monge-Amp\`ere equation, Math.
Res. Lett. 19 (2012), no. 4, 939-946.

\noindent
[Y] S.T. Yau, {\sl On the Ricci curvature of a compact K\"ahler manifold and the complex Monge-Amp\`ere
equation. I}, Comm. Pure Appl. Math. {\bf  31} (1978), 339-411.

\end{document}